%% file: main.tex
\documentclass[11pt]{amsart}
\usepackage{amssymb,graphics}
\usepackage[colorlinks]{hyperref}
\usepackage[all]{xy}

\title{Hodge Theory and Derived Categories of Cubic Fourfolds}
\author[N.~Addington]{Nicolas Addington}
\email{adding@math.duke.edu}
\address{\noindent Department of Mathematics \\
Duke University, Box 90320 \\
Durham, NC 27708-0320 \\
United States}
\author[R.~P.~Thomas]{Richard Thomas}
\email{richard.thomas@imperial.ac.uk}
\address{Department of Mathematics \\
Imperial College London \\
London SW7 2AZ \\
United Kingdom}

\newcommand \A {\mathcal A}
\newcommand\B{\mathcal B}
\newcommand \C {\mathbb C}
\newcommand \cC {\mathcal C}

\newcommand \I {\mathcal I}
\newcommand \K {\mathcal K}

\newcommand \M {\mathcal M}

\newcommand \OO {\mathcal O}
\renewcommand \P {\mathbb P}
\newcommand \p {\mathfrak p}
\newcommand \Q {\mathbb Q}
\newcommand \R {\mathbb R}
\newcommand \cS {\mathcal S}
\newcommand \X {\mathcal X}
\newcommand \Z {\mathbb Z}

\DeclareMathOperator \im {im}
\renewcommand \phi \varphi
\newcommand \Ktop {K_{\rm top}}
\newcommand \Kalg {K_{\rm alg}}
\newcommand \Knum {K_{\rm num}}
\DeclareMathOperator \hcf {hcf\,}

\DeclareMathOperator \disc {disc}
\DeclareMathOperator \sig {sig}
\DeclareMathOperator \Bl {Bl}

\DeclareMathOperator \rank {rank}
\DeclareMathOperator \ch {ch}   
\DeclareMathOperator \td {td}
\let \Re \undefined
\DeclareMathOperator \Re {Re}

\DeclareMathOperator \pr {pr}
\DeclareMathOperator \Hom {Hom}
\DeclareMathOperator \Ext {Ext}
\newcommand \mar {{\rm mar}}
\newcommand \lev {{\rm lev}}

\DeclareMathOperator \Gal {Gal}
\newcommand \prim {{\rm prim}}

\newcommand \HH {H\!H}
\newcommand \hook {\mathbin{\lrcorner}}
\DeclareMathOperator \Spec {Spec}
\DeclareMathOperator \At {At}
\DeclareMathOperator \sHom {\mathcal Hom}

\newtheorem*{thm*}{Theorem}
\newtheorem{thm}{Theorem}
\numberwithin{thm}{section}
\newtheorem{prop}[thm]{Proposition}
\newtheorem{lem}[thm]{Lemma}
\newtheorem{cor}[thm]{Corollary}
\theoremstyle{definition}
\newtheorem{defn}[thm]{Definition}
\newtheorem*{rmk*}{Remark}

\numberwithin{equation}{section}

\begin{document}

\input intro
\input mukai-lattice
\input kuz2hast
\input non-empty
\input open
\input richard
\input hodge_conj

\newcommand \httpurl [1] {\href{http://#1}{\nolinkurl{#1}}}
\bibliographystyle{plain}
\bibliography{cubic_4folds}
\vspace{5mm}
\end{document}

%% file: intro.tex

\begin{abstract}
Cubic fourfolds behave in many ways like K3 surfaces. Certain cubics -- conjecturally, the ones that are rational -- have specific K3 surfaces associated to them geometrically. Hassett has studied cubics with K3 surfaces associated to them at the level of Hodge theory, and Kuznetsov has studied cubics with K3 surfaces associated to them at the level of derived categories.

These two notions of having an associated K3 surface should coincide. We prove that they coincide generically: Hassett's cubics form a countable union of irreducible Noether-Lefschetz divisors in moduli space, and we show that Kuznetsov's cubics are a dense subset of these, forming a non-empty, Zariski open subset in each divisor.

\end{abstract}

\maketitle

\section{Introduction}
It has long been noted that there are remarkable similarities between cubic fourfolds and K3 surfaces \cite{BD, Rap, voisin_thesis}. Let $X$ denote a smooth complex cubic hypersurface in $\P^5$. Both the Hodge structure and the derived category of $X$ decompose into some trivial pieces and a K3-like piece.

\subsection*{Hodge theory} The Hodge-theoretic viewpoint is due to Hassett \cite{hassett_thesis}.  The Hodge diamond of $X$ is
\begin{equation*}
\begin{array}{ccccc} 
& & 1 \\ & & 1 \\ 0 & 1 & 21 & 1 & 0 \\ & & 1 \\ & & \phantom{.}1.
\end{array}
\end{equation*}
Removing the powers $h^i$ of the hyperplane class, we are left with the primitive cohomology:
\[ \begin{array}{ccccc}
0 & 1 & 20 & 1 & 0.
\end{array} \]
This looks like (a Tate twist of) $H^2$ of a K3 surface $S$, but the intersection forms have different signatures: $(20,2)$ for $H^4_\prim(X,\Z)$ versus $(19,3)$ for $H^2(S,\Z)(-1)$.  However, one can often find a codimension-1 sub-Hodge structure of signature $(19,2)$ common to both. For the K3 surface it is $H^2_\prim(S,\Z)(-1)$, the orthogonal to some primitive ample class $\ell$. For the cubic $X$ it is $\langle h^2,T\rangle^\perp$, the subspace of $H^4_\prim(X,\Z)$ orthogonal to an integral $(2,2)$-class $T\in H^{2,2}(X,\Z)$.

Up to automorphisms of the lattice, both of these situations are governed by a single positive integer $d$. For the K3 surface it is the degree $d=\ell^2$ of the ample class $\ell$, while for the cubic fourfold it is the discriminant $d=\disc\langle h^2,T\rangle$. (Without loss of generality we always choose $T$ so that $\langle h^2,T\rangle$ is a \emph{primitive} sublattice of $H^4(X,\Z)$.) By \cite[Thm.~1.0.1]{hassett_thesis}, the cubics possessing such a $T\in H^{2,2}(X,\Z)$ form an irreducible divisor $\cC_d$ in the 20-dimensional moduli space $\cC$ of cubics, non-empty if and only if
\begin{equation} \label{even} \tag{$*$}
d > 6 \mathrm{\ and\ } d \equiv 0 \mathrm{\ or\ } 2 \pmod 6.
\end{equation}
Moreover, by \cite[Thm.~5.1.3]{hassett_thesis} there exists a polarised K3 surface $(S,\ell)$ with
\begin{equation} \label{hass}
H^2_\prim(S,\Z)(-1)\ \cong\ \langle h^2,T\rangle^\perp
\end{equation}
if and only if $d$ satisfies the further condition\footnote{For $d$ even (which is implied by \eqref{even}), this strange-looking numerical condition turns out to be equivalent to $d$ being the norm of a primitive vector in the $A_2$ lattice $\left( \begin{smallmatrix} 2 & -1 \\ -1 & 2 \end{smallmatrix} \right)$.}
\begin{equation} \label{numerical_condition} \tag{$**$}
\text{$d$ is not divisible by 4, 9, or any odd prime $p \equiv 2 \pmod 3$.}
\end{equation}

\subsection*{Derived category}
The derived category viewpoint is due to Kuznetsov \cite{kuznetsov}.  The line bundles $\OO_X$, $\OO_X(1)$, and $\OO_X(2)$ form an exceptional collection in $D(X)$, and the right orthogonal
\begin{eqnarray*} \A_X &:=& \langle \OO_X, \OO_X(1), \OO_X(2) \rangle^\perp \\ &:=& \{E\in D(X)\colon R\Hom(\OO_X(i),E)=0\ \mathrm{for}\ i=0,1,2\}
\end{eqnarray*}
looks like the derived category of a K3 surface in that it has the same Serre functor and Hochschild (co)homology. In general, $\A_X$ should be thought of as a non-commutative K3 surface: it is a deformation of the derived category of a genuine K3 surface but, for general $X$, we will see it lacks point-like objects -- that is, objects $E \in\A_X$ with $\Ext^*_{\A_X}
(E,E) \cong \Ext^*_\text{K3}(\OO_\text{point}, \OO_\text{point})$. We will call $\A_X$ \emph{geometric} if $\A_X \cong D(S)$ for some projective\footnote{``Projective'' may be redundant: $\A_X$ is \emph{saturated} in the sense of \cite{bvdb}, and it is expected that if $S$ is not projective then $D(S)$ is not saturated \cite[Rem.~5.6.2]{bvdb}.}  K3 surface $S$.

\subsection*{Rationality}
Although it is not strictly relevant to our paper, a strong motivation is the question of which cubic fourfolds are rational.

If $X$ is rational then by ideas of Clemens-Griffiths the transcendental cohomology of $X$ must come from the cohomology of a surface $S$.  Given the shape of $H^*(X)$ it is natural to expect $S$ to be a K3 surface.  Conversely, Hassett showed that many of those $X$ whose transcendental cohomology is that of a K3 surface are rational. So Harris and Hassett asked if having an associated K3 surface in the Hodge-theoretic sense should be equivalent to rationality. (They were cautious about conjecturing this, however. In fact no cubic fourfold has yet been shown to be irrational.)

For similar reasons, Kuznetsov \cite{kuznetsov} conjectured that $X$ is rational if and only if $X$ has an associated K3 surface in the derived category sense: i.e.\ $\A_X$ is geometric.  He showed that this is true of the known rational cubics.  His conjecture has attracted a great deal of interest recently \cite{ms2, bfk, abbv}.

\subsection{Our results}
At the very least, Hassett and Kuznetsov's conditions should be the
same. (This would give a non-commutative extension of the derived Torelli theorem of Mukai-Orlov, which states that two projective K3 surfaces have equivalent derived categories if and only if they have isomorphic Mukai lattices.) We prove this in the generic case:

\begin{thm} \label{MAIN}
If $\A_X$ is geometric then $X \in \cC_d$ for some $d$ satisfying \eqref{numerical_condition}.  Conversely, for each $d$ satisfying \eqref{numerical_condition}, the set of cubics $X \in \cC_d$ for which $\A_X$ is geometric is a Zariski open dense subset.
\end{thm}
\noindent Of course we fully expect that this open dense subset is all of $\cC_d$. But%
, as we discuss in Section \ref{end},%
\ the issue of taking the ``right" limit of a family of objects of the derived category makes proving this difficult at present.

There are no such problems with taking limits of algebraic cycles, however. Applied to the Mukai vector of the Fourier-Mukai kernel of Theorem \ref{MAIN}, this leads to the following result in Section \ref{hodge_conj}, generalising
\eqref{hass}.

\begin{thm} \label{algcycle} 
Fix a cubic fourfold $X\in\cC_d$ for some $d$ satisfying \eqref{numerical_condition}. Then there exists a polarised K3 surface $S$ of degree $d$ and an algebraic cycle in $A^3(S\times X)_\Q$ which induces a Hodge isometry $H^2_\prim(S,\Z)(-1)\to\langle h^2, T \rangle^\perp$.
\end{thm}

We also strengthen this as follows.

\begin{thm} \label{HODGE}
Let $S$ be a projective K3 surface and $X$ a cubic fourfold. If a Hodge class $Z \in H^{3,3}(S \times X, \Q)$ induces a Hodge isometry $T(S)(-1)\xrightarrow\sim T(X)$ between \emph{integral} transcendental lattices then $Z$ is algebraic.
\end{thm}

\subsection{Strategy} \label{plan}
An outline of the paper is as follows.
In Section \ref{mukai-lattice} we define a Mukai lattice for the category $\A_X$, that is, a weight-2 Hodge structure which is isomorphic to the usual Mukai lattice $H^*(S,\Z)$ when $\A_X \cong D(S)$. We do this using topological K-theory; in fact most of the lattice-theoretic complications in this subject stem from the difference between the integral structures on cohomology and K-theory.

We spend some time relating this Mukai lattice to $H^4(X,\Z)$. In Section \ref{kuz2hast} this enables us to interpret \eqref{numerical_condition} in much more natural K-theoretic terms. Namely, $X \in \cC_d$ for some $d$ satisfying
\eqref{numerical_condition} if and only if there are classes $\kappa_1,\kappa_2
\in\Kalg(\A_X)$ that behave like the classes of a skyscraper sheaf and ideal sheaf of a point on a K3 surface: $\chi(\kappa_1,\kappa_1)=\chi(\kappa_2,\kappa_2)=0$ and $\chi(\kappa_1,\kappa_2)=1$.

Morally, this is the reason for the equivalence of Hassett's and Kuznetsov's conditions. We should now ``just" produce the K3 surface $S$ as a moduli space of objects in $\A_X$ of class $\kappa_1\in\Ktop(\A_X)$. To do this we first work in $\cC_d\cap\cC_8$, then use deformation theory to reach a Zariski open subset of $\cC_d$.

The advantage of $\cC_8$ is that we have Kuznetsov's description of $\A_X$ for $X\in\cC_8$ as the derived category of \emph{twisted} sheaves on a K3 surface $S$. (Notice that this does not make $\A_X$ geometric unless the twisting cocycle vanishes, and indeed $d=8$ does not satisfy \eqref{numerical_condition}.)  So in Section \ref{non-empty} we first use some lattice theory to show that $\cC_d\cap\cC_8$ contains cubics $X$ whose extra discriminant-$d$ class forces the twisting cocyle to vanish, so $\A_X\cong D(S)$ is indeed geometric. However, this equivalence is the ``wrong" one, expressing $S$ as a moduli space of objects in a class different from $\kappa_1$, so it does not deform out of $\cC_8$.

But we have still gained something: since $\A_X$ is geometric we can now work in $D(S)$, where the powerful results of Mukai \cite{mukai_tata} give us a moduli space of stable sheaves on $S$ whose class in $\Ktop(\A_X)$ is $\kappa_1$. Replacing $S$ by this moduli space in Section \ref{open} gives us the ``correct" equivalence $\A_X\cong D(S)$.

In Section \ref{2} we deform $X$ inside $\cC_d$. There is a corresponding deformation of $S$, using the Torelli theorem and the close relationship between the cohomologies of $X$ and $S$. That their Hodge structures remain related (in the sense of Hassett) in the family means that the cohomological obstruction to deforming the equivalence -- or, more precisely, its Fourier-Mukai kernel -- vanishes. That is, the Mukai vector of the kernel remains of type $(p,p)$ in the family. Using some Hochschild (co)homology theory set up in Section \ref{1} we show that this implies the vanishing of the obstructions to deforming the kernel to any order. By algebraicity, the equivalence deforms over a Zariski open subset of $\cC_d$.

\subsection*{Acknowledgements}
We thank K.~Buzzard, T.~Gee, H.~C.~Johannsson, J.~Ludwig, A.~Skorobogatov and F.~Thorne for generous assistance with quadratic forms, and B.~Hassett, D.~Huybrechts and C.~Voisin for expert advice on this project. We also thank T.~Coates, J.~Francis and C.~Westerland for help with topology, W.~Donovan, E.~Macr\`i, E.~Markman, S.~Mehrotra, E.~Segal and P.~Stellari for useful discussions and correspondence, two thorough referees who suggested many improvements, and A.~Bayer for a correction to \S\ref{finally}. This work was supported by an EPSRC Programme Grant no.\ EP/G06170X/1.

\subsection*{Notation} All our varieties are smooth and complex projective. We usually use $X$ to denote a cubic fourfold, and $S$ a K3 surface. Both have torsion-free cohomology, so it makes sense to let $H^{p,p}(X,\Z)$ denote the intersection $H^{p,p}(X)\cap H^{2p}(X,\Z)$ inside $H^{2p}(X,\C)$.

The bounded derived category of coherent sheaves on $Y$ is denoted $D(Y)$. In Section \ref{2.2}, where $Y$ might denote a smooth family over an Artinian base, $D(Y)$ denotes the bounded derived category of perfect complexes.

%% file: mukai-lattice.tex

\section{A Mukai lattice for \texorpdfstring{$\A_X$}{AX}} \label{mukai-lattice}

In this section we introduce a weight-2 Hodge structure which we call the \emph{Mukai lattice} of $\A_X$, and relate it to the lattice $H^4(X,\Z)$.

\subsection{The Mukai vector and topological K-theory}
Fix smooth projective varieties $Y$ and $Z$. The Mukai vector of an object $E\in D(Y)$ is
\begin{equation} \label{Mukv}
v(E)=\ch(E)\td(Y)^{1/2}\ \in\, H^*(Y,\Q).
\end{equation}
Any object $P\in D(Y\times Z)$ induces a Fourier-Mukai functor $\Phi_P\colon D(Y)\to D(Z)$ defined by the formula
\[ \Phi_P(\ \cdot\ ):=\pi_{Z*}(\pi_Y^*(\ \cdot\ )\otimes P), \]
where all functors are derived unless otherwise stated. There is an induced map $\Phi_P^H$ on rational cohomology given by
\begin{equation} \label{HP}
\Phi_P^H(\ \cdot\ ):=\pi_{Z*}\big(\pi_Y^*(\ \cdot\ )\cup v(P)\big) \colon\ H^*(Y,\Q) \longrightarrow H^*(Z,\Q).
\end{equation}
If $Y$ and $Z$ are K3 surfaces then $\Phi_P^H$ acts on \emph{integral} cohomology. But in general it is a different integral structure on rational cohomology that $\Phi_P^H$ respects: \emph{topological K-theory}, or rather its image under the Mukai vector.  Since this is rarely used in algebraic geometry we give a brief review, using \cite{ah1} as our main reference. We take
\[ \Ktop(Y) = \Ktop^0(Y) \oplus \Ktop^1(Y), \]
although for the spaces we consider $\Ktop^1$ vanishes.  The Mukai vector \eqref{Mukv} induces an isomorphism of vector spaces $\Ktop(Y) \otimes \Q \to H^*(Y,\Q)$. A map $f\colon Y \to Z$ induces, in addition to the usual pullback map $f^*\colon \Ktop(Z) \to \Ktop(Y)$, a pushforward map $f_*\colon \Ktop(Y) \to \Ktop(Z)$ satisfying a projection formula and a Grothendieck-Riemann-Roch formula.  By \cite{ah2} it is compatible with the pushforward on $\Kalg$.  Applying this to the projection $Y\times Z\to Z$ we find that any object $P \in D(Y \times Z)$ induces a map $\Phi^K_P\colon \Ktop(Y) \to \Ktop(Z)$ compatible with the usual induced maps:
\[ \xymatrix{
D(Y) \ar[d]_{\Phi_P} \ar[r] 
& \Kalg(Y) \ar[d]_{\Phi^K_P} \ar[r] 
& \Ktop(Y) \ar[d]_{\Phi^K_P} \ar[r]^(.46)v 
& H^*(Y,\Q) \ar[d]_{\Phi^H_P} \\
D(Z) \ar[r] & \Kalg(Z) \ar[r] & \Ktop(Z) \ar[r]^(.46)v & H^*(Z,\Q).\!
} \]
Using the map $p\colon Y\to\mathrm{point}$ we can define an Euler pairing on $\Ktop(Y)$,
\begin{equation} \label{pair}
\chi(\kappa_1,\kappa_2) := p_*(\kappa_1^\vee \otimes \kappa_2)\in
\Ktop(\mathrm{point})\cong\Z,
\end{equation}
which agrees with the usual Euler pairing on $\Kalg$ and satisfies
\begin{equation} \label{serre_duality}
\chi(\kappa_1, \kappa_2) = (-1)^{\dim Y} \chi(\kappa_2, \kappa_1 \otimes\omega_Y),
\end{equation}
as can be seen from Riemann-Roch on $\Ktop(Y)$. If $P_L, P_R \in D^b(Y \times Z)$ are the kernels inducing the left and right adjoints of $\Phi_P$ then $\Phi_{P_L}^K$ and $\Phi_{P_R}^K$ are left and right adjoint to $\Phi_P^K$ under the Euler pairing.

Since K3 surfaces and cubic fourfolds have torsion-free cohomology, the following will be very useful to us.
\begin{thm}[Atiyah-Hirzebruch {\cite[\S2.5]{ah1}}] \label{ahss}
If $H^*(Y,\Z)$ is torsion-free then
\begin{enumerate}
\item $\Ktop(Y)$ is torsion-free and $v\colon \Ktop(Y) \to H^*(Y,\Q)$ is injective.
\item For any $\kappa \in \Ktop(Y)$, the leading term of $v(\kappa)$ is integral: writing $v(\kappa) = v_j + v_{j+1} + \dotsb$ with $v_i \in H^i(Y,\Q)$, we have $v_j \in H^j(Y,\Z)$.
\item For any $j$ and any $v_j \in H^j(Y,\Z)$ there is a $\kappa \in \Ktop(Y)$ with $v(\kappa) = v_j + \text{higher-degree terms}$.
\end{enumerate}
\end{thm}

\subsection{The Mukai lattice of a K3 surface}
If $S$ is a K3 surface then $v\colon \Ktop(S) \to H^*(S,\Z)$ is an isomorphism, so we can rephrase the usual Mukai lattice in terms of $\Ktop$. That is, the Mukai lattice of $S$ is the Abelian group $\Ktop(S)$ endowed with the Euler pairing \eqref{pair} and the weight-2 Hodge structure $\Ktop(S)\otimes\C = \bigoplus_{p+q=2} \widetilde H^{p,q}(S)$, where
\begin{align*}
\widetilde H^{2,0}(S) &:= v^{-1}(H^{2,0}(S)), \\
\widetilde H^{1,1}(S) &:= v^{-1}(H^{0,0}(S) \oplus H^{1,1}(S) \oplus H^{2,2}(S)), \\
\widetilde H^{0,2}(S) &:= v^{-1}(H^{0,2}(S)).
\end{align*}
Note that Mukai originally used the opposite sign for the pairing.

Since $S$ satisfies the integral Hodge conjecture and $\Ktop(S)$ is torsion-free, we can identify the numerical K-theory $\Knum(S) := \Kalg(S)/\ker\chi$ with
\[ \im\big(\Kalg(S) \to \Ktop(X)\big)\,=\ \Ktop(S) \cap \widetilde H^{1,1}(S). \]

\subsection{A Mukai lattice for \texorpdfstring{$\A_X$}{AX}} \label{KAX}
The above suggests the following definition.
\begin{defn}
If $X$ is a cubic fourfold, the \emph{Mukai lattice of $\A_X$} is the Abelian group
\[ \Ktop(\A_X) := \{ \kappa \in \Ktop(X) : \chi([\OO_X(i)], \kappa) = 0 \text{ for } i = 0, 1, 2 \} \]
endowed with the Euler pairing \eqref{pair} and the weight-2 Hodge structure $\Ktop(\A_X)\otimes\C = \bigoplus_{p+q=2} \widetilde H^{p,q}(\A_X)$, where
\begin{align*}
\widetilde H^{2,0}(\A_X) &:= v^{-1}(H^{3,1}(X)), \\
\widetilde H^{1,1}(\A_X) &:= v^{-1}(H^{0,0}(X) \oplus H^{1,1}(X) \oplus H^{2,2}(X) \oplus H^{3,3}(X) \oplus H^{4,4}(X)), \\
\widetilde H^{0,2}(\A_X) &:= v^{-1}(H^{1,3}(X)).
\end{align*}
\end{defn}
Since $X$ satisfies the Hodge conjecture \cite{zucker} and $\Ktop(X)$ is torsion-free, we can identify the numerical K-theory $\Knum(\A_X) := \Kalg(\A_X)/\ker\chi$ with
\[ \im\big(\Kalg(\A_X) \to \Ktop(\A_X)\big)\,\subseteq\ \Ktop(\A_X) \cap \widetilde H^{1,1}(\A_X). \]
In Proposition \ref{lemmar} we will see that this inclusion is an equality, essentially because $X$ satisfies the integral Hodge conjecture \cite{voisin_Z_hodge}.

If $S$ is a K3 surface and $\Phi\colon D(S) \to D(X)$ is a fully faithful functor with image $\A_X$ we get a Hodge isometry
\[ \Phi^K\colon \Ktop(S) \xrightarrow\sim\Ktop(\A_X) \]
with inverse induced by the right adjoint of $\Phi$. Since there are cubics $X$ with $\A_X \cong D(S)$ (and all smooth cubics are deformation equivalent), the Euler pairing on $\Ktop(\A_X)$ is symmetric, which was not obvious \emph{a priori} \eqref{serre_duality}. It is abstractly isomorphic to the even unimodular lattice $U^4 \oplus {E_8}^2$.

\subsection{Relation to \texorpdfstring{$H^4(X,\Z)$}{H4(X,Z)}}
The numerical K-theory $\Knum(\A_X)$ always contains at least two linearly independent classes\footnote{Of course for $\A_X$ to be geometric we need at least \emph{three} linearly independent classes in $\Knum(\A_X)$, corresponding to the classes in $\Knum(S)$ of $\OO_\text{point}$, $\I_\text{point}$ and a polarisation.}, given by projecting  $\OO_\text{line}$ and $\OO_\text{point}$ into $\A_X$. In fact we will find it more convenient to use $\OO_\text{line}(1)$ and $\OO_\text{line}(2)$ instead. Precisely, let
\begin{eqnarray*}
\pr_i\colon \Ktop(X) &\!\!\longrightarrow\!\!& \Ktop(X), \\
\kappa &\mapsto& \kappa - \chi\big([\OO_X(i)],\,\kappa\big) \cdot [\OO_X(i)],
\end{eqnarray*}
be the projection onto $[\OO_X(i)]^\perp$, so
\[ \pr := \pr_0 \circ \pr_1 \circ \pr_2 \]
projects $\Ktop(X)$ onto $\Ktop(\A_X)$. Now define
\[
\lambda_1 := \pr\,[\OO_\text{line}(1)] \text{ and } \lambda_2 := \pr\,[\OO_\text{line}(2)].
\]
Calculating their Euler pairing we find that they generate the sublattice
\begin{equation} \label{A_2}
{-A_2} = \begin{pmatrix} -2 & 1 \\ 1 & -2 \end{pmatrix}\ \subset\ \Knum(\A_X)\ \subset\ \Ktop(\A_X),
 \end{equation}
which is negative definite, and primitive since $\disc(-A_2) = 3$ is square-free.  The orthogonal to these classes in $\Ktop(\A_X)$ is comparable to $H^4_\prim(X,\Z)$: in the first case we have removed $[\OO_X(i)]$, $i=0,1,\dotsc,4$, from $\Ktop(X)$ while in the second we have removed $h^i$, $i=0,1,\dotsc,4$, from $H^*(X,\Z)$.  Remarkably, they are in fact integrally isometric via the Mukai vector.

\begin{prop} \label{orthogonal}
The Mukai vector $v\colon \Ktop(\A_X) \to H^*(X,\Q)$ takes \linebreak $\langle \lambda_1,\lambda_2 \rangle^\perp \subset \Ktop(\A_X)$ isometrically onto $\langle h^2 \rangle^\perp \subset H^4(X,\Z):$
\begin{equation} \label{perpy}
\langle \lambda_1,\lambda_2 \rangle^\perp\ \cong\ H^4_\prim(X,\Z).
\end{equation}
More generally, if $\kappa_1, \dotsc, \kappa_n \in \Ktop(\A_X)$ then the Mukai vector takes \linebreak $\langle \lambda_1,\lambda_2,\kappa_1,\dotsc,\kappa_n \rangle^\perp$ isometrically onto $\langle h^2, c_2(\kappa_1), \dotsc, c_2(\kappa_n) \rangle^\perp$.
\end{prop}
\begin{proof}
Since the cohomology of $X$ is so simple,
\[ H^*(X,\Q)\cong\langle1,h,h^2,h^3,h^4\rangle\oplus H^4_\prim(X,\Q), \]
and since the Todd class of $X$ is a linear combination of the $h^i$, we find that for any $\kappa \in \Ktop(X)$,
\begin{align*}
&\ \ \kappa\in\langle\lambda_1,\lambda_2\rangle^\perp\subset\Ktop(\A_X) \\
&\iff \kappa\in\langle\OO_X, \OO_X(1), \OO_X(2), \lambda_1,\lambda_2\rangle^\perp\subset\Ktop(X) \\
&\iff \kappa\in\langle\OO_X, \OO_X(1), \OO_X(2), \OO_X(3), \OO_X(4)\rangle^\perp\subset\Ktop(X) \\
&\iff v(\kappa)\in\langle1,h,h^2,h^3,h^4\rangle^\perp\subset H^*(X,\Q) \\
&\iff v(\kappa)\in\langle h^2\rangle^\perp\subset H^4(X,\Q) \\
&\iff v(\kappa)\in\langle h^2\rangle^\perp\subset H^4(X,\Z),
\end{align*}
where the last line comes from the leading term of $v(\kappa)$ being integral.

So we must show that any $T\in\langle h^2\rangle^\perp\subset H^4(X,\Z)$ is the image of some $\tau\in\langle\lambda_1,\lambda_2\rangle^\perp\subset\Ktop(\A_X)$. By Theorem \ref{ahss} there is a class $\tau'\in \Ktop(X)$ with $v(\tau') = T + \text{higher-degree terms}$. Since $T.h^2 = 0$ we have
\[\chi(\OO_X(t), \tau') = at + b\quad\text{for some }\ a, b \in \Z,\]
by Riemann-Roch on $\Ktop$.  Set $\tau:= \tau' + a[\OO_\text{line}]-(a+b)[\OO_\text{point}]$. Then $\chi(\OO_X(t), \tau) = 0$ for all $t$, so $\tau \in \langle \lambda_1,\lambda_2 \rangle^\perp \subset \Ktop(\A_X)$.

For the more general statement, note that if $v(\tau)\in\langle h^2\rangle^\perp$ then $v(\tau)$ equals its leading term $-c_2(\tau)$, so another application of Riemann-Roch gives 
\[
\chi(\tau, \kappa_i) = c_2(\tau) \cdot c_2(\kappa_i) \quad
\text{for any } \tau \in \langle \lambda_1,\lambda_2 \rangle^\perp \subset \Ktop(\A_X). \qedhere
\]
\end{proof}

For a unimodular lattice $L$ and a non-degenerate sublattice $M$, the obvious inclusion $M^\perp\subset L/M$ is not an isomorphism unless $M$ is unimodular. Thus
$\langle\lambda_1,\lambda_2\rangle^\perp \subsetneq \Ktop(\A_X)/\langle \lambda_1,\lambda_2 \rangle$ and $\langle h^2 \rangle^\perp \subsetneq H^4(X,\Z)/\langle h^2\rangle$. But the isomorphism \eqref{perpy} of Proposition \ref{orthogonal} extends to the bigger groups too:

\begin{prop} \label{lemmar} The second Chern class $c_2$ descends to an isomorphism of free Abelian groups
\[ \bar c_2\colon \frac{\Ktop(\A_X)}{\langle\lambda_1,\lambda_2\rangle}
\xymatrix{\ar[r]^\sim &} \frac{H^4(X,\Z)}{\langle h^2\rangle}\,. \]
The preimage of $H^{2,2}(X,\Z)/\langle h^2 \rangle$ is the image of $\Kalg(\A_X)$.  In particular, $\Ktop(\A_X) \cap \widetilde H^{1,1}(\A_X) = \Knum(\A_X)$.
\end{prop}

\begin{proof}
Taking $c_2$ and projecting gives a map $\Ktop(\A_X) \to H^4(X,\Z)/\langle h^2\rangle$. This is a group homomorphism because
\[
c_2(\kappa_1 + \kappa_2) = c_2(\kappa_1) + c_1(\kappa_1)c_1(\kappa_2) + c_2(\kappa_2) \]
and the middle term lies in $\langle h^2\rangle$.

Clearly $\langle\lambda_1,\lambda_2\rangle$ lies in the kernel. Conversely, take $\kappa \in \Ktop(\A_X)$ with $c_2(\kappa)$ a multiple of $h^2$. Pick two distinct smooth hyperplane sections $H_1,H_2$ of $X$. The leading terms of the Mukai vectors of
\begin{equation} \label{lisst}
[\OO_X], [\OO_{H_1}(1)], [\OO_{H_1 \cap H_2}(2)], [\OO_\text{line}(1)], [\OO_\text{point}] \end{equation}
are the $\Z$-generators of $\langle 1,h,h^2,h^3,h^4\rangle_\Q\cap H^4(X,\Z)$:
\[ 1,\ h,\ h^2,\ h^3/3,\ h^4/3 \]
respectively. Therefore by an induction on the degree of the leading term of $v(\kappa)$ -- which is always integral by Theorem \ref{ahss} -- we find that $\kappa$ is in the integral span of \eqref{lisst}.  Applying $\pr$ to the classes \eqref{lisst} gives
\[ 0,\ 0,\ 0,\ \lambda_1,\ \lambda_2-\lambda_1, \]
so $\kappa=\pr(\kappa)$ is a linear combination of $\lambda_1,\lambda_2$.  Thus $\bar c_2$ is injective.

To see that $\bar c_2$ is surjective, let $T \in H^4(X,\Z)$.  By Theorem \ref{ahss}, there is a $\tau \in \Ktop(X)$ with $v(\tau) = -T +{}$higher-degree terms. In particular, $c_2(\tau)=T$.  Then $\pr(\tau) \in \Ktop(\A_X)$ differs from $\tau$ by a linear combination of $[\OO_X]$, $[\OO_X(1)]$, and $[\OO_X(2)]$ whose Chern classes are all multiples of $h^i$.  Thus $c_2(\pr(\tau))$ differs from $c_2(\tau) = T$ by a multiple of $h^2$. \smallskip

If $T \in H^{2,2}(X,\Z)$ we can take $\tau$ above to be the image of an algebraic class, since Voisin has proved the \emph{integral} Hodge conjecture for cubic fourfolds \cite{voisin_Z_hodge}. Thus the natural map $\Kalg(\A_X)\to\Ktop(\A_X)\cap \widetilde H^{1,1}(\A_X)$ is surjective.
\end{proof}

\begin{rmk*}
The generic cubic 4-fold has $H^{2,2}(X,\Z) = \langle h^2 \rangle$, so by
Proposition \ref{lemmar}, $\Knum(\A_X) = \langle \lambda_1, \lambda_2 \rangle \cong -A_2$ is negative definite.  Thus $\A_X$ cannot be equivalent to $D(S)$, or even to a derived category of twisted sheaves $D(S,\alpha)$, because $[\OO_\text{point}]$ is isotropic in $\Knum(S,\alpha)$. In Section \ref{kuz2hast} we will get more precise information along these lines.
\end{rmk*}

\begin{prop} \label{two_lattices}
Given $\kappa_1, \dotsc, \kappa_n \in \Ktop(\A_X)$, consider the sublattices
\begin{gather*}
M_H:=\langle h^2,c_2(\kappa_1),\dotsc,c_2(\kappa_n)\rangle\subset H^4(X,\Z),\\
M_K:=\langle\lambda_1,\lambda_2,\kappa_1,\dotsc,\kappa_n \rangle \subset \Ktop(\A_X).
\end{gather*}
\begin{enumerate}
\setlength \itemindent {-\leftmargin}
\addtolength \itemindent {2em}
\item A class $\kappa \in \Ktop(\A_X)$ is in $M_K$ if and only if $c_2(\kappa) \in M_H$,
\item $M_H$ is primitive if and only $M_K$ is,
\item $M_H$ is non-degenerate if and only if $M_K$ is, and
\item if $M_K \subset \Knum(\A_X)$ then $M_K$ and $M_H$ are non-degenerate.
\end{enumerate}
Moreover, when $M_H$ and $M_K$ are non-degenerate,
\begin{enumerate}
\setlength \itemindent {-\leftmargin}
\addtolength \itemindent {2em}
\setcounter{enumi}{4}
\item $\sig M_H=(r,s)$ if and only if $\sig M_K=(r-1,s+2)$, and
\item $M_H$ and $M_K$ have the same discriminant.
\end{enumerate}
\end{prop}

\begin{proof}
For the first two claims we use the isomorphism $\bar c_2$ of Proposition \ref{lemmar}. Under the two projections
\[
M_K\subset \Ktop(\A_X) \xymatrix{\ar[r]^(.47){\pi_K}&}
\frac{\Ktop(\A_X)}{\langle\lambda_1,\lambda_2\rangle}\ \cong\ 
\frac{H^4(X,\Z)}{\langle h^2\rangle} \xymatrix{&\ar[l]_(.43){\pi_H}}
H^4(X,\Z)\supset M_H
\]
the sublattices $M_K$ and $M_H$ project to the same subgroup
\[
\overline M:=\langle\bar\kappa_1,\dotsc,\bar\kappa_n\rangle\ \cong\ 
\langle \bar c_2(\kappa_1),\dotsc,\bar c_2(\kappa_n) \rangle.
\]
Moreover $M_K$ and $M_H$ contain the kernels of the projections, so in fact
\begin{equation} \label{projn}
M_K=\pi_K^{-1}(\overline M)\quad\text{and}\quad M_H=\pi_H^{-1}(\overline M),
\end{equation}
from which it follows that for any $\kappa\in\Ktop(\A_X)$,
\[ \kappa\in M_K\iff\bar c_2(\kappa)\in\overline M\iff c_2(\kappa)\in M_H. \]
Similarly \eqref{projn} implies that
\[ M_K\text{ is primitive }\iff\overline M\text{ is primitive }\iff M_H\text{ is primitive}, \]
and more generally
\begin{equation} \label{i}
i(M_K)=i(\overline M)=i(M_H),
\end{equation}
where $i(\ \cdot\ )$ denotes the index of each sublattice in its saturation.

The remaining claims follow from the isomorphism $M_K^{\perp}\cong M_H^\perp$ of Proposition \ref{orthogonal}, as follows.

For claim (3), use $\rank(M_H) = \rank(\overline M) + 1$ and $\rank(M_K) = \rank(\overline M)+2$ from \eqref{projn}, and that $M_H$ (respectively $M_K$) is non-degenerate if and only if $\rank(M_H^\perp) = 23 - \rank(M_H)$ (respectively $\rank(M_K^\perp) = 24 - \rank(M_K)$). Claim (4) then follows from (3) and the Hodge-Riemann bilinear relations, which imply that $M_H\subset H^{2,2}(X,\Z)$ is positive definite.

For claim (5), observe that $H^4(X,\Z)$ has signature $(21,2)$ by the Hodge-Riemann bilinear relations, and $\Ktop(\A_X)\cong U^4\oplus {E_8}^2$ has signature $(20,4)$, so if $M_H$ and $M_K$ are non-degenerate then
\[ \sig M_H\!=\!(r,s)\!\!\iff\!\!\sig(M_H^\perp\!=\!M_K^\perp)\!=\!(21-r,2-s)\!\!\iff\!\!
\sig M_K\!=\!(r-1,s+2). \]
In particular, the discriminants of $M_H$ and $M_K$ have the same sign.

For claim (6) we use the fact \cite[Cor.~1.6.2]{nikulin} that if $M$ is a primitive sublattice of a unimodular lattice then $\lvert\disc(M)\rvert = \lvert\disc(M^\perp)\rvert$.  Both $\Ktop(\A_X) \cong U^4 \oplus {E_8}^2$ and $H^4(X,\Z)$ are unimodular; the latter by Poincar\'e duality.  Therefore, letting $i$ denote any of the indices \eqref{i}, we have
\begin{align*} 
\lvert \disc(M_K) \rvert
= i^2\lvert \disc(M_K^{\perp\perp}) \rvert
= i^2\lvert \disc&(M_K^\perp = M_H^\perp) \rvert \\
&= i^2\lvert \disc(M_H^{\perp\perp}) \rvert
= \lvert \disc(M_H) \rvert. \qedhere
\end{align*}
\end{proof}

%% file: kuz2hast.tex

\section{Interpretation of the numerical condition} \label{kuz2hast}

If $\A_X \cong D(S)$ for some K3 surface $S$ then the classes $\kappa_1=[\OO_\text{point}]$ and $\kappa_2=[\I_\text{point}]$ span a copy of $U= \left(\begin{smallmatrix} 0 & 1 \\ 1 & 0 \end{smallmatrix}\right)$ in $\Knum(\A_X)$. Combined with the following result this proves the easier direction of Theorem \ref{MAIN}.
\begin{thm} \label{day}
Let $X$ be a cubic fourfold.  The following are equivalent:
\begin{enumerate}
\item $X \in \cC_d$ for some $d$ satisfying \eqref{numerical_condition},
\item $\Knum(\A_X)$ contains a copy of the hyperbolic plane $U = \left(\begin{smallmatrix} 0 & 1 \\ 1 & 0 \end{smallmatrix}\right)$.\end{enumerate}
\end{thm}
\noindent Kuznetsov \cite[Prop.~4.8]{kuznetsov} proved a special case of this, namely that a generic $X \in \cC_8$ does not have a copy of $U$ in $\Knum(\A_X)$.

\begin{proof}[Proof of Theorem \ref{day}] 
By definition, $X \in \cC_d$ if and only if there is a primitive sublattice $M_H \subset H^{2,2}(X,\Z)$ of rank 2 and discriminant $d$ that contains $h^2$.  Thus by Propositions \ref{lemmar} and \ref{two_lattices} we see that $(1)$ is equivalent to
\begin{enumerate}
\item[$(1')$] There is a primitive sublattice $M_K \subset \Knum(\A_X)$ of rank 3 and discriminant $d$ that contains $\langle\lambda_1,\lambda_2\rangle$, such that $d$ satisfies \eqref{numerical_condition}.
\end{enumerate}
Hassett \cite[Prop.~5.1.4]{hassett_thesis} has shown that $d = \disc(M_H)$ satisfies \eqref{numerical_condition} if and only if there is a primitive embedding of the rank-1 lattice $(-d)$ into $U^3 \oplus {E_8}^2$ such that the orthogonal is isomorphic to ${M_H}^\perp \subset H^4(X,\Z)$.  So by Proposition \ref{orthogonal} we see that $(1')$ is equivalent to
\begin{enumerate}
\item[$(1'')$] There is a primitive sublattice $M_K \subset \Knum(\A_X)$ of rank 3 and discriminant $d$ that contains $\langle\lambda_1,\lambda_2\rangle$, such that ${M_K}^\perp \subset \Ktop(\A_X)$ is isomorphic to $(-d)^\perp$ for some primitive embedding
\begin{equation} \label{hassett_embedding}
(-d) \subset U^3 \oplus {E_8}^2.
\end{equation}
\end{enumerate}
Now we argue as follows:

\smallskip
$(1'') \Rightarrow (2)$.  Adding $U$ to both sides of \eqref{hassett_embedding}, we have primitive embeddings of $M_K$ and $U \oplus (-d)$ into $\Ktop(\A_X)\cong U^4 \oplus {E_8}^2$ with isomorphic orthogonals.  Therefore $M_K$ and $U \oplus (-d)$ have the same discriminant quadratic form \cite[Cor.~1.6.2]{nikulin}.  Since one contains $U$, they are isomorphic \cite[Cor.~1.13.4]{nikulin}.  Thus we have $U \subset M_K \subset \Knum(\A_X)$.

\smallskip
$(2) \Rightarrow (1'')$ or  $(1')$.  Conversely, suppose that $\kappa_1, \kappa_2 \in \Knum(\A_X)$ span a copy of $U$.  Since $\langle \lambda_1, \lambda_2 \rangle \cong -A_2$ is negative definite and $U$ is indefinite we see that the rank of $\langle \lambda_1, \lambda_2, \kappa_1, \kappa_2 \rangle$ is 3 or 4.  We prove $(1'')$ in the first case and $(1')$ in the second.
\smallskip

{\bf Rank 3.}  In this case let $M_K$ be the saturation of $\langle \lambda_1, \lambda_2, \kappa_1, \kappa_2 \rangle$, and let $d = \disc(M_K)$.  We have inclusions
\[ U \subset M_K \subset \Ktop(\A_X) \cong U^4 \oplus {E_8}^2. \]
Since $U$ is unimodular, $M_K$ splits as the direct sum of $U$ and the orthogonal to $U$ in $M_K$, so the latter is a rank-1 lattice of discriminant $-d$.  Similarly the orthogonal to $U$ in $\Ktop(\A_X)$ is an even unimodular lattice of signature $(19,3)$, hence is isomorphic to $U^3 \oplus {E_8}^2$.  Thus ${M_K}^\perp \subset \Ktop(\A_X)$ is isomorphic to $(-d)^\perp \subset U^3 \oplus {E_8}^2$, which gives $(1'')$.
\smallskip

{\bf Rank 4.}  In this case our argument is unfortunately rather longer and more \emph{ad hoc}.  We will show that there are integers $x$ and $y$, not both zero, such that the discriminant of
\begin{equation} \label{sublattice_for_qf}
\langle \lambda_1, \lambda_2, x \kappa_1 + y \kappa_2 \rangle
\end{equation}
satisfies \eqref{numerical_condition}.  To see that this implies $(1')$, let $M_K$ be the saturation of \eqref{sublattice_for_qf}; then $\disc(M_K)$ divides the discriminant of \eqref{sublattice_for_qf}, so if the latter satisfies \eqref{numerical_condition} then the former does too.  Note that \eqref{sublattice_for_qf} has rank 3 unless $x=y=0$.

We consider the quadratic form
\[ Q(x,y) := \begin{cases}
\disc \langle \lambda_1, \lambda_2, x \kappa_1 + y \kappa_2 \rangle & \text{if $x\ne0$ or $y\ne0$,} \\
0 & \text{if }x=y=0.
\end{cases} \]
This is positive definite: since $c_2(x \kappa_1 + y \kappa_2) \in H^{2,2}(X)$, Proposition \ref{two_lattices} and the Hodge-Riemann bilinear relations give $Q(x,y) > 0$ unless $x=y=0$.

Writing the Euler pairing on $\langle \lambda_1, \lambda_2, \kappa_1, \kappa_2 \rangle$ as
\begin{equation} \label{euler_pairing_rank_4}
\begin{pmatrix}
-2 & 1 & k & m \\
1 & -2 & l & n \\
k & l & 0 & 1 \\
m & n & 1 & 0
\end{pmatrix},
\end{equation}
we find that
\[ Q(x,y) = A x^2 + B xy + C y^2, \]
where the coefficients $A$, $B$, and $C$ are the following minors of \eqref{euler_pairing_rank_4}:
\[ A = \begin{vmatrix}
-2 & 1 & k \\
1 & -2 & l \\
k & l & 0
\end{vmatrix},
\qquad
B = 2 \begin{vmatrix}
-2 & 1 & m \\
1 & -2 & n \\
k & l & 1
\end{vmatrix},
\qquad
C = \begin{vmatrix}
-2 & 1 & m \\
1 & -2 & n \\
m & n & 0
\end{vmatrix}. \]
Let $h = \hcf(A,B,C)$ be the highest common factor of these coefficients, let
\[ a = A/h,\quad b = B/h,\quad c = C/h, \]
and
\[ q(x,y) = ax^2 + bxy + cy^2. \]
To finish the proof it is enough to show $h$ satisfies \eqref{numerical_condition} and that $q$ represents a prime $p \equiv 1 \pmod 3$, i.e.\ there are $x,y\in\Z$ such that $q(x,y)=p$.  These are proved in the next two results.
\end{proof}

\begin{lem}
$\hcf(A,B,C)$ is even and satisfies \eqref{numerical_condition}.
\end{lem}
\begin{proof}
Since the lattice $-A_2$ appears in the top left hand corner \eqref{euler_pairing_rank_4}, it is convenient to phrase things in terms of the \emph{Eisenstein integers} $\Z[\omega]\subset\C$, where $\omega = e^{2 \pi i/3} = \tfrac{-1 + \sqrt{-3}}2$. Endowing $\Z[\omega]$ with the integral bilinear form $(\alpha,\gamma)\mapsto 2\Re(\alpha \bar\gamma)$ gives a lattice isomorphic to $A_2$.  So we replace the variables $k,l,m,n$ of \eqref{euler_pairing_rank_4} with the Eisenstein integers
\[ \alpha:= k - l\omega, \qquad \gamma:= m - n\omega. \]
In these variables we find that
\[ A = 2 |\alpha|^2,
\qquad B = 4 \Re(\alpha \bar\gamma) + 6,
\qquad C = 2 |\gamma|^2, \]
where $|\alpha|^2$ denotes the usual Euclidean norm.  In particular $h=\hcf(A,B,C)$ is even. We will use the following standard facts about the ring $\Z[\omega]$: it is a principal ideal domain in which
\begin{enumerate}
\item $\sqrt{-3} = 1 + 2\omega$ is prime,
\item every prime $p \in \Z$ with $p \equiv 2 \pmod 3$ is prime in $\Z[\omega]$, and
\item\label{norm_is_not_2} every $\beta \in \Z[\omega]$ has $|\beta|^2 \equiv 0$ or $1 \pmod 3$.
\end{enumerate}
The last of these is easily checked by listing the elements of $\Z[\omega]$ modulo 3. \medskip

Suppose that $p$ is an odd prime with $p \equiv 2 \pmod 3$. Then if $p\mid A$ we have
\[ p \mid 2 \alpha \bar\alpha \Longrightarrow p \mid \alpha \Longrightarrow p \mid 2 \Re(\alpha \bar\gamma) \Longrightarrow p \nmid B = 4 \Re(\alpha \bar\gamma) + 6 \Longrightarrow p \nmid h, \]
as required.  Similarly we find that $4 \nmid h$.

Suppose that $9=(\sqrt{-3})^4$ divides both $A = 2 \alpha \bar\alpha$ and $C = 2 \gamma \bar\gamma$. Then $\sqrt{-3}$ divides $\alpha$ and $\gamma$ twice, so
\[ 9 \mid 2 \Re(\alpha \bar\gamma)\ \Longrightarrow\ 9 \nmid B = 4 \Re(\alpha \bar\gamma) + 6\ \Longrightarrow\ 9 \nmid h. \qedhere \]
\end{proof}

\begin{prop} \label{primitive_qf} The primitive positive definite quadratic form $q(x,y) = ax^2 + bxy + cy^2$ represents a prime $p \equiv 1 \pmod 3$.
\end{prop}
\begin{proof}
Define $D := b^2 - 4ac$.  There are two cases: $3 \mid D$ and $3 \nmid D$.

\medskip
{\bf Case 1: $\bf 3 \mid D$.}  By simply listing the quadratic forms modulo 3 for which $D(q)\equiv0\pmod3$, we find that $q$ represents only
\begin{itemize}
\item $0 \pmod 3$ \hspace{13mm} (when $q\equiv0\pmod 3$), or
\item 0 and $1 \pmod 3$ \quad ($q\equiv x^2,\ y^2,\ x^2+xy+y^2,\ x^2+2xy+y^2$), or \item 0 and $2 \pmod 3$ \quad ($q\equiv2x^2,\ 2y^2,\ 2x^2+2xy+2y^2,\ 2x^2+xy+2y^2$),
\end{itemize}
but not all of $0,1,2\pmod3$. The first case occurs only when $q$ is not primitive, so we ignore it. In the second case, since any primitive positive definite form represents a prime (in fact infinitely many \cite[Thm.~9.12]{cox}) it must represent a prime $p \equiv 1 \pmod 3$, as required.

So it is enough to show that $a \equiv 1 \pmod 3$ or $c \equiv 1 \pmod 3$ since both imply that $q$ represents $1\pmod3$, putting us in the second case. \smallskip

Suppose first that $|\alpha|^2,|\gamma|^2$ are both $0\pmod3$. Then
\[ \sqrt{-3}\mid\alpha,\gamma\ \Longrightarrow\ 3\mid B\ \Longrightarrow\ 3\mid h. \]
Since $h$ is even and satisfies \eqref{numerical_condition} we see that the integer $h/6\equiv1\pmod3$.
Therefore, writing $\alpha = \sqrt{-3} \alpha'$ and $\gamma = \sqrt{-3} \gamma'$, we have
\[ |\alpha'|^2 = a(h/6)\equiv a\pmod3 \quad\mathrm{and}\quad |\gamma'|^2=c(h/6)
\equiv c\pmod3. \]
So it is enough to show that one of $|\alpha'|^2,|\gamma'|^2$ is $1 \pmod 3$.  But if not then by fact \eqref{norm_is_not_2} above we have $\sqrt{-3} \mid \alpha'$ and $\sqrt{-3} \mid \gamma'$, so $3 \mid 2 \Re(\alpha' \bar\gamma')$, so
\[ D\equiv D(h/6)^2 = (2 \Re(\alpha' \bar\gamma') + 1)^2 - 4|\alpha'|^2|\gamma'|^2 \equiv 1 \pmod 3, \]
contradicting our assumption that $3 \mid D$. \smallskip

So by fact \eqref{norm_is_not_2} we are left with the case where at least one of $|\alpha|^2,|\gamma|^2$ is $1\pmod3$. Then $h$ is not divisible by 3, so $h/2\equiv1\pmod3$ by \eqref{numerical_condition}. Hence
\[ |\alpha|^2 = a(h/2)\equiv a\pmod3,\quad |\gamma|^2 = c(h/2)\equiv c\pmod3, \]
and at least one of $a,\,c$ is $1 \pmod 3$.\medskip

{\bf Case 2: $\bf 3 \nmid D$.} In this case we adapt the usual proof that $q$ represents infinitely many primes to show that it represents infinitely many congruent to $1 \pmod 3$.  We follow \cite[Thm.~9.12]{cox}.  Consider $K := \Q(\sqrt D)$, the order $\OO \subset \OO_K$ of discriminant $D$, and its ring class field $L$.  Via Artin reciprocity, $q$ corresponds to an element $\sigma_0 \in \Gal(L/K)$.  Let $\langle \sigma_0 \rangle$ denote the conjugacy class of its image in $\Gal(L/\Q)$.  To show that $q$ represents infinitely many primes it is enough to show that the Dirichlet density of
\begin{equation} \label{dirichlet_thing}
\left\{ \text{$p$ prime} : \text{$p$ is unramified in $L$}, \left(\frac{L/\Q}{p}\right) = \langle \sigma_0 \rangle \right\}
\end{equation}
is positive, which is true by \v Cebotarev's density theorem.  To show that $q$ represents infinitely many primes congruent to $1 \pmod 3$, we adapt this as follows.

First we claim that 3 is unramified in $L$. Certainly 3 is unramified in $K$, since $3 \nmid D$.  Either $3 \OO_K$ is prime in $\OO_K$, or $3 \OO_K = \p \bar\p$ for some prime $\p$ of $\OO_K$.  In the first case it is enough to show that $3 \OO_K$ is unramified in $L$, and in the second that $\p$ and $\bar\p$ are unramified in $L$.

By \cite[p.~180]{cox}, if a prime of $\OO_K$ is ramified in $L$ then it divides $f \OO_K$, where $f = [ \OO_K : \OO ]$ is the conductor of $\OO$. Therefore its norm divides $N(f \OO_K) = f^2$, hence divides $f^2 d_K = D$, where $d_K$ is the discriminant of $K$.  In the two cases above we have $N(3 \OO_K) = 9$ and $N(\p) = N(\bar\p) = 3$, neither of which divides $D$.

Since 3 is unramified in $L$ we have $\sqrt{-3} \notin L$, so let $L' = L(\sqrt{-3})$.  Then
\[ \Gal(L'/\Q) \cong \Gal(L/\Q) \times \Gal(\Q(\sqrt{-3})/\Q) = \Gal(L/\Q) \times \{ \pm 1 \}. \]
Saying that $p \equiv 1 \pmod 3$ is equivalent to saying that the Artin symbol $\left(\frac{\Q(\sqrt{-3})/\Q}{p}\right) = 1$.  Thus we replace \eqref{dirichlet_thing} with
\[ \left\{ \text{$p$ prime} : \text{$p$ is unramified in $L'$}, \left(\frac{L'/\Q}{p}\right) = \langle \sigma_0 \rangle \times 1 \right\} \]
whose Dirichlet density is again positive by \v Cebotarev's density theorem.
\end{proof}

%% file: non-empty.tex

\section{Non-emptiness} \label{non-empty}

\begin{thm} \label{41}
Suppose that $\cC_d$ is non-empty, i.e.\ that $d$ satisfies \eqref{even}.  Then there is a $X \in \cC_d \cap \cC_8$ such that $\A_X$ is geometric.
\end{thm}

A cubic $X \in \cC_8$ contains a plane $P$ by \cite[\S3]{voisin_thesis}. The linear system of hyperplanes containing $P$ defines a map $\Bl_P(X) \to \P^2$ whose fibres are quadric surfaces $Q$ (a point of $\P^2$ -- the dual of the linear system -- corresponds to two hyperplanes whose intersection is the reducible cubic surface $Q\cup P$). We also let $Q = h^2-P$ denote the class in $X$ (rather than $\Bl_P(X)$) of a smooth fibre. We have $P^2=3$ and $Q^2=4$.

The fibres degenerate over a sextic curve $C \subset \P^2$.  If it is smooth then the double cover of $\P^2$ branched over $C$ is a K3 surface $S$, and there is a natural Brauer class $\alpha \in \operatorname{Br}(S)$.  Kuznetsov \cite[\S4]{kuznetsov} has shown that $\A_X \cong D(S, \alpha)$, the derived category of $\alpha$-twisted sheaves on $S$.

To prove the theorem it is enough to find such an $X\in\cC_8$ with a class $T\in H^{2,2}(X,\Z)$ such that the sublattice $\langle h^2, T \rangle$ is primitive of discriminant $d$ and $T.Q=1$: by \cite[Prop.~4.7]{kuznetsov} this implies that $\alpha=0$, so $\A_X$ is geometric.

So we fix $L:=H^4(Y,\Z)$ for some cubic fourfold $Y$ containing a plane, with corresponding classes $h^2, P\in L$, and $Q = h^2 - P$. We will first produce a suitable $T\in L$, then use Laza and Looijenga's description of the image of the period map to produce another cubic $X$ such that these classes are in $H^{2,2}(X,\Z)$ .  We begin with the following.

\begin{lem} \label{hast_prop}
For any $n>0$ with $n \equiv 5 \pmod 8$, there is a $T \in L$ such that $T.Q=1$ and $\langle h^2,Q,T\rangle\subset L$ is primitive of discriminant $n$ with intersection pairing
\begin{equation} \label{16k-3}
\qquad \begin{pmatrix}
3 & 2 & 0 \\
2 & 4 & 1 \\
0 & 1 & 2k
\end{pmatrix} \qquad\mathrm{when\ }n = 16k-3,\ \mathrm{or}
\end{equation}
\begin{equation} \label{16k+5}
\qquad \begin{pmatrix}
3 & 2 & 1 \\
2 & 4 & 1 \\
1 & 1 & 2k+1
\end{pmatrix} \quad\mathrm{when\ }n = 16k+5.
\end{equation}
In particular, $\langle h^2,T \rangle \subset L$ is primitive of discriminant $6k$ in \eqref{16k-3} and $6k+2$ in \eqref{16k+5}.
\end{lem}

\begin{proof}
Hassett {\cite[Prop.~4.1 and proof of Lemma~4.4]{hassett_rational}} already gives us a $T$ such that $\langle h^2,Q,T\rangle\subset L$ is primitive of discriminant $n$ and $T.Q = 1$. The intersection pairing on $\langle h^2,Q,T \rangle$ is of the form
\[ \begin{pmatrix}
3 & 2 & 4a+b \\
2 & 4 & 1 \\
4a+b & 1 & *
\end{pmatrix} \]
for some $0 \le b \le 3$.  Replacing $T$ with $T - 2a h^2 + aQ$, this becomes
\begin{equation} \label{lvs}
\begin{pmatrix}
3 & 2 & b \\
2 & 4 & 1 \\
b & 1 & *
\end{pmatrix}.
\end{equation}
Since changing $T$ to $h^2 - T$ replaces $b$ with $3-b$ in \eqref{lvs}, we may assume without loss of generality that $b$ is either 0 or 1.

Now $\disc\langle h^2,T\rangle = 3T^2 - b^2$ is necessarily even; one way to see this is to observe that $T':=3T - bh^2$ lies in $\langle h^2 \rangle^\perp$, so $T'^2 = 9T^2 - 3b^2$ is necessarily even by \cite[Prop.~2.1.2]{hassett_thesis} or Proposition \ref{orthogonal}.  Thus if $b=0$ then $T^2$ is even, so setting $2k = T^2$ we get \eqref{16k-3}. If $b=1$ then $T^2$ is odd, so setting $2k+1 = T^2$ we get \eqref{16k+5}.  
\end{proof}

\begin{proof}[Proof of Theorem \ref{41}.]
Since the case $d=8$ is known we fix $d > 8$ with $d \equiv 0$ or $2 \pmod 6$. We apply Lemma \ref{hast_prop} with $n = 16k-3$ (if $d = 6k$ with $k \ge 2$) or $n = 16k+5$ (if $d = 6k+2$ with $k \ge 2$) to produce a $T\in L$ such that the intersection form on $\langle h^2,Q,T\rangle$ is \eqref{16k-3} or \eqref{16k+5}. Thus for any $T' = x h^2 + y Q + z T$ we have
\begin{equation} \label{rqfs}
\disc\langle h^2, T' \rangle = 8y^2 + 6yz + 6k z^2 \quad\mathrm{or}\quad
8y^2 + 2yz + (6k+2) z^2
\end{equation}
respectively. By the theory of reduced quadratic forms \cite[\S5.2]{baker}, the smallest non-zero value taken by these forms is 8, at only $y = \pm 1, z = 0$. In particular,
\begin{equation} \label{nexists}
\nexists\ T'\in\langle h^2,Q,T\rangle\ \mathrm{with}\ \disc\langle h^2, T'\rangle = 2\ \mathrm{or}\ 6.
\end{equation}

\medskip
Next we produce a $\sigma \in L\otimes\C$ with
\[ \sigma^2 = 0,\ \sigma\bar\sigma < 0\ \ \mathrm{and}\ \ L\cap\sigma^\perp = \langle h^2,Q,T \rangle. \]
Recall that the signature of $L$ is $(21,2)$.  Choose a 21-dimensional positive definite real subspace $V\subset L\otimes\R$ such that $V\cap L = \langle h^2, Q, T \rangle$.  Then the pairing is negative definite on $V^\perp \subset L\otimes\R$, so any non-zero $\sigma \in V^\perp\otimes\C \subset L\otimes\C$ has $\sigma \bar\sigma < 0$.  Moreover, the line $\P(V^\perp\otimes\C) \subset \P(L\otimes\C)$ necessarily meets the quadric $\{ \sigma^2 = 0 \}$.

Then by Laza and Looijenga's description of the image of the period map \cite[Thm.~1.1]{laza}, \cite[Thm.~3.1]{looijenga}, the condition \eqref{nexists} guarantees the existence of a cubic $X$ and an isomorphism $L\cong H^4(X,\Z)$ preserving $h^2$ and taking $\sigma$ to a generator of $H^{3,1}(X)$.  In particular, $\langle h^2, Q, T \rangle=H^{2,2}(X,\Z)$ by the construction of $V=\langle\sigma,\bar\sigma\rangle^\perp$, so $X \in \cC_d$.  Moreover the cohomology class $P = h^2-Q$ corresponds to an actual plane $P\subset X$ by \cite[\S3]{voisin_thesis}.

Finally we claim that the discriminant sextic of $\Bl_P(X) \to \P^2$ is smooth. Voisin \cite[\S1]{voisin_thesis} observes that it is enough to show that there is no other plane $P' \subset X$ meeting $P$, and that distinct planes have independent homology classes. Such a $P'\in\langle h^2, Q, T \rangle$ would have $\disc\langle h^2,P'\rangle=8$, which by our discussion of reduced quadratic forms following \eqref{rqfs} implies that $P' = x h^2 \pm Q$.  From $h^2.P' = 1$ we deduce that $P' = h^2 - Q = P$, as desired.
\end{proof}

%% file: open.tex

\section{Setup for deformation theory} \label{open}

\subsection{Modification of the equivalence}
We now modify the equivalence $\Phi\colon D(S) \xrightarrow\sim \A_X$ of Theorem \ref{41} to the ``right one" when $d$ satisfies \eqref{numerical_condition}.

We want the equivalence to deform as $X$ deforms through $\cC_d$. A necessary condition -- which we shall find in Section \ref{2} is also sufficient -- is that $c_2(\Phi(\OO_\text{point}))$ and $c_2(\Phi(\I_\text{point}))$ should lie in $\langle h^2, T \rangle$ so that they remain algebraic throughout $\cC_d$.

Since $v(\OO_\text{point}), v(\I_\text{point}) \in H^*(S,\Z)$ are orthogonal to $H^2_\prim(S,\Z)$ and $\Phi^H$ respects the Mukai pairing on cohomology, it is equivalent to ask that $\Phi^H$ takes $H^2_\prim(S,\Z)$ to $\langle h^2,T \rangle^\perp$. This is what Proposition \ref{mukai-orlov} below achieves. \medskip

Let $X_0 \in \cC_d$ for some $d$ satisfying \eqref{numerical_condition}. Then by Hassett \eqref{hass} there is a polarised K3 surface $(S_0,\ell_0)$ and a Hodge isometry
\[ \phi\colon H^2_\prim(S_0,\Z)(-1)\longrightarrow\langle h^2,T \rangle^\perp. \]

\begin{prop} \label{mukai-orlov}
If $\A_{X_0}$ is geometric then we can choose an equivalence $\Phi_{P_0}\colon D(S_0) \xrightarrow\sim \A_{X_0}$ so that $\Phi_{P_0}^H\colon H^*(S_0,\Q) \to H^*(X_0,\Q)$ extends $\phi$.
\end{prop}

\begin{proof}
First we extend $\phi$ to a Hodge isometry $\Ktop(S_0) \to \Ktop(\A_{X_0})$.   Recall that $\Ktop(S_0)$ and $\Ktop(\A_{X_0})$ are both isomorphic as lattices to $U^4 \oplus {E_8}^2$, and that $-H^2_\prim$ and $\langle h^2, T \rangle^\perp$ are naturally primitive sublattices, the latter by Proposition \ref{orthogonal}.  The orthogonal $(-H^2_\prim)^\perp \subset \Ktop(S_0)$ contains a copy of $U$, so by \cite[Thm.~1.14.4]{nikulin} any two primitive embeddings $-H^2_\prim \hookrightarrow U^4 \oplus {E_8}^2$ differ by an automorphism of $U^4 \oplus {E_8}^2$.  Thus the lattice isomorphism $\phi\colon {-H^2_\prim} \to \langle h^2, T \rangle^\perp$ extends to a lattice isomorphism $\tilde\phi\colon \Ktop(S_0) \to \Ktop(\A_{X_0})$.  Since $\phi$ takes $H^{2,0}(S_0)$ to $H^{3,1}(X_0)$ we see that the extension $\tilde\phi$ does as well, so $\tilde\phi$ is a Hodge isometry.

Now since $\A_{X_0}$ is geometric there is a K3 surface $S$ and a Fourier-Mukai equivalence $D(S) \xrightarrow\sim \A_{X_0}$.  Let $\psi \colon \Ktop(S) \to \Ktop(\A_{X_0})$ be the induced Hodge isometry.  If $\psi^{-1} \circ \tilde\phi \colon \Ktop(S_0) \to \Ktop(S)$ does not preserve the natural orientation of the four negative directions, change the sign of $\tilde\phi$ on the copy of $U \subset \Ktop(S_0)$ spanned by $[\OO_\text{point}]$ and $[I_\text{point}]$; since $U$ has one positive and one negative direction, $\psi^{-1} \circ \tilde\phi$ now preserves the orientation of the negative directions.  The new $\tilde\phi$ is still an extension of $\phi$.

Now the strong form of the derived Torelli theorem \cite[Cor.~10.13]{huybrechts_fm} gives an equivalence $D(S_0) \xrightarrow\sim D(S)$ such that the induced map on $\Ktop$ is $\psi^{-1} \circ \tilde\phi$.  Composing with $D(S) \xrightarrow\sim \A_{X_0}$ gives the required $\Phi_{P_0}\colon D(S_0) \xrightarrow\sim \A_{X_0}$.
\end{proof}

\subsection{Construction of the families} \label{construct}
As described in Section \ref{plan}, we now want to deform $X_0$ inside $\cC_d$ and take a corresponding deformation of $S_0$.  Hassett \cite[Cor.~5.2.4]{hassett_thesis} has shown that the isometry $\phi\colon H^2_\prim(S_0,\Z)(-1) \xrightarrow\sim \langle h^2,T \rangle^\perp$ yields an open immersion of the moduli space of cubics with marked $\langle h^2, T \rangle$ into the moduli space of polarised K3 surfaces.  But these are not fine moduli spaces, which is a technical inconvenience for us in Section \ref{2.2}, so we pass to a finite cover.

\begin{prop} \label{state}
Let $X_0$, $d$, $S_0$, $\ell_0$, and $P_0$ be as in Proposition \ref{mukai-orlov}.  There is a smooth quasi-projective variety $\cC_d^\lev$ with the following properties. There is a family of cubics
\[ p_X\colon \X \longrightarrow \cC_d^\lev, \]
such that the classifying map $\pi\colon \cC_d^\lev \to \cC_d$ is a finite cover, and a family
\[ p_S\colon (\cS, \ell) \longrightarrow \cC_d^\lev \]
of smooth polarised K3 surfaces. There is a point $0 \in \cC_d^\lev$ mapping to $X_0 \in \cC_d$, over which the fibre of $(\cS, \ell)$ is $(S_0, \ell_0)$. Finally, the embedding $\Phi_{P_0}^H\colon H^*(S_0, \Q) \hookrightarrow H^*(X_0,\Q)$ extends to an embedding of local systems
\[ \xymatrix{ R p_{S*} \Q\ \ar@{^(->}[r] & \,R p_{X*} \Q }\]
whose complexification takes the period point $H^{2,0}(S_t) \subset H^*(S_t,\C)$ to the period point $H^{3,1}(X_t) \subset H^*(X_t,\C)$ for all $t \in \cC_d^\lev$.
\end{prop}

\begin{proof}
Consider the lattices
\begin{align*}
L_0:=\langle h^2 \rangle^\perp &\subset H^4(X_0,\Z)=:L, \\
\Lambda_0 := \langle \ell_0 \rangle^\perp &\subset H^2(S_0,\Z)=:\Lambda,
\end{align*}
and the associated local period domains
\begin{align*}
D&:= \{ \sigma \in \P(L_0 \otimes \C) : \sigma^2 = 0,\,\sigma\bar\sigma < 0 \}, \\
\Delta&:= \{ \sigma \in \P(\Lambda_0 \otimes \C) : \sigma^2 = 0,\,\sigma\bar\sigma > 0 \}.
\end{align*}

Let $\cC^\mar$ be the moduli space of marked cubic fourfolds, that is, cubics $X$ with an isomorphism $H^4(X,\Z) \cong L$ preserving $h^2$.  The period map embeds $\cC^\mar$ as an open subset of $D$ by \cite{voisin_thesis}.  Similarly, let $\K_d^\mar$ be the moduli space of marked polarised K3 surfaces of degree $d$; the period map embeds $\K_d^\mar$ as an open subset of $\Delta$.

We identify both moduli spaces $\cC^\mar,\,\K_d^\mar$ with their images in the period domains.  Both carry universal families: for $\K_d^\mar$ the Torelli theorem allows one to glue local universal deformations, as is detailed in \cite[Exp.~XIII]{mep} or
\cite[\S6.3.3]{huybrechts_K3}, and the same argument works for $\cC^\mar$ thanks to Voisin's Torelli theorem of \cite{voisin_thesis}. (While \cite{voisin_thesis} does not emphasise that every Hodge isometry $H^4(X_1,\Z) \to H^4(X_2,\Z)$ preserving $h^2$ is induced by a \emph{unique} isomorphism $X_1 \to X_2$, it does enough to prove it.)

Next define\footnote{Note that this is not what Hassett calls $\cC_d^\mar$ in \cite{hassett_thesis}; his is a quotient of ours by the group $G$ below.} $\cC_d^\mar = \cC^\mar \cap T^\perp$, which parametrises marked cubics with a distinguished integral $(2,2)$-class of discriminant $d$.  Then the lattice isomorphism $\phi\colon {-\Lambda_0} \to \langle h^2, T \rangle^\perp$ induces an isomorphism of period domains $\phi\colon \Delta \to D \cap T^\perp$.  Moreover we claim that
\begin{equation} \label{last}
\phi^{-1}(\cC_d^\mar) \subset \K_d^\mar.
\end{equation}
Indeed, $\K_d^\mar$ is the complement in $\Delta$ of the hyperplanes $\delta^\perp$, where $\delta \in \Delta$ is any class with $\delta^2 = -2$.  But a smooth cubic has no integral $(2,2)$ class $\delta$ with $h^2.\delta = 0$ and $\delta^2 = 2$ by \cite[\S4,~Prop.~1]{voisin_thesis}, or by noting that $\disc\langle h^2,\delta \rangle$ would be 6 which contradicts \eqref{even}.

Thus over $\cC_d^\mar$ we get both a family $\X_d^\mar$ of cubics (by restricting the universal family over $\cC^\mar$) \emph{and} a family $\cS^\mar$ of K3 surfaces. By construction their period points coincide under $\phi$ and hence under $\Phi_{P_0}^H$.

Since $\cC_d^\mar$ is not quasi-projective, we take an appropriate quotient.  Consider the subgroup $G$ of $O(L)$ fixing $h^2$ and $T$.  By \cite[Cor.~1.5.2]{nikulin} this can be identified with the subgroup of $O(\langle h^2, T \rangle^\perp)$ acting trivially on the discriminant group of $\langle h^2, T \rangle^\perp$; this is sometimes called the stable orthogonal group.  Via $\phi$ this is identified with subgroup of $O(\Lambda_0)$ acting trivially on the discriminant group of $\Lambda_0$, and hence with the subgroup of $O(\Lambda)$ fixing $\ell_0$.  Observe that the embedding \eqref{last} is $G$-equivariant, and that the action of $G$ on $\cC_d^\mar$ lifts to an action on the families $\X_d^\mar$ and $\cS^\mar$ by the strong form of the global Torelli theorems.

The action of $G$ on $\cC_d^\mar$ has finite stabilisers: the stabiliser of a point in $\cC_d^\mar$ is the automorphism group of the corresponding cubic, which is finite \cite{mm}. Therefore we fix $N \in \Z$ large enough that the subgroup
\[ G_N = \{ g \in G : g \equiv 1\!\! \pmod N \} \]
is torsion-free, and set $\cC_d^\lev = \cC_d^\mar/G_N$.  Now $G_N$ acts freely on $\cC_d^\mar$, so $\cC_d^\lev$ is smooth, and the families descend to $\cC_d^\lev$: take $\X = \X_d^\mar/G_N$ and $\cS = \cS^\mar/G_N$.  Since $G_N$ is a finite-index subgroup of $O(\Lambda_0)$, the quotient $\cC_d^\lev$ is quasi-projective by results of Baily and Borel.

Letting $0 \in \cC_d^\lev$ be the projection of $X_0$ (or its period point) in $\cC_d^\mar$ completes the proof.
\end{proof}

%% file: richard.tex

\section{(Hochschild) (co)homology} \label{1}

The deformation theory of $D(Y)$ is governed by the Hochschild cohomology $\HH^*(Y)$, which for $Y=S$ a K3 surface is isomorphic, as an ungraded vector space, to both its Hochschild homology $\HH_*(S)$ and its de Rham cohomology $H^*(S)$. This means we can often reduce deformation/obstruction theory questions to Hodge theory. We are guided by the deformation theory for Fourier-Mukai equivalences in the wonderful papers \cite{Toda, HMS}\footnote{There is also related work \cite{abp} using deformation quantisation.} but with a number of simplifications and one generalisation (to Fourier-Mukai functors which are only fully faithful). We also use $T^1$-lifting methods to simplify the higher order deformation theory. We therefore give a self-contained account.

\subsection*{Notation: families}
Although we suppress mention of it throughout, we work relative to a complex affine base $B$ in this section and the next. So for instance a ``smooth projective variety" means a smooth projective family over $B$. All products, diagonals, etc.\ are taken relative to $B$; thus $\Delta_S$ denotes the diagonal in $S\times_BS$, etc.  Usually $B$ can be taken to be $\Spec\C$, but in Section \ref{2.2} we will apply our results to the case $B=\Spec A$ for an Artinian local $\C$-algebra $A$. The notation is unaffected; we simply replace our scalars $\C$ by $A=\OO_B$.

Since we sometimes work with smooth families over \emph{singular} bases $B$, it is important that by $D(Y)$ we will always mean the bounded derived category of \emph{perfect} complexes of coherent sheaves on $Y$. (Note that the structure sheaf of the diagonal $\OO_{\Delta_Y}$ is indeed perfect in $D(Y\times_BY)$.) For $B$ Artinian all cohomologies and Exts will therefore be finite dimensional.

\subsection{Hochschild cohomology}
Suppose that $S$ and $X$ are smooth projective varieties and that $P\in D(S\times X)$ is the Fourier-Mukai kernel of a \emph{fully faithful} embedding $D(S)\hookrightarrow D(X)$.

The usual Fourier-Mukai convolution product $*$, corresponding to composition of Fourier-Mukai functors, gives the functors
\begin{equation} \label{diagram}
\xymatrix@R=5pt{
D(S\times S) \ar[r]^{P*}& D(S\times X) & D(X\times X), \ar[l]_{\ *P} \\
\qquad \OO_{\Delta_S}\quad \ar@{|->}[r]& \qquad P\qquad & \quad\OO_{\Delta_X}.\qquad \ar@{|->}[l].}
\end{equation}
The first induces a map
\begin{equation} \label{SP}
\Ext^*(\OO_{\Delta_S},\OO_{\Delta_S}) \xrightarrow{\ P*} \Ext^*(P,P)
\end{equation}
which is an isomorphism. In fact\footnote{Alternatively use the right adjoint $\Phi_R$ \eqref{Radjo}, which is a left inverse to the fully faithful $\Phi_P$.
Thus $R * P\cong\OO_{\Delta_X}$ and $R*$ is also right adjoint to $P*$. Therefore \eqref{SP} is the isomorphism
$\Ext^*(\OO_{\Delta_S},\OO_{\Delta_S})\cong\Ext^*(\OO_{\Delta_S}, R * P)\cong\Ext^*(P,P)$.}
\[ \pi_{1*} \sHom_{S\times S}(\OO_{\Delta_S},\OO_{\Delta_S}) \xrightarrow{\ P*}  \pi_{1*} \sHom_{S\times X}(P,P) \]
is a quasi-isomorphism because the full and faithful condition ensures that it is a quasi-isomorphism when restricted to any $s\in S$, where it is the map $\Ext^*(\OO_s,\OO_s)\to\Ext^*(P_s,P_s)$.

The right hand map of \eqref{diagram} gives a map $\Ext^*(\OO_{\Delta_X},\OO_{\Delta_X})
\to\Ext^*(P,P)$. Combining with the inverse of the isomorphism \eqref{SP} gives
\begin{equation} \label{HH*}
\HH^*(X):=\Ext^*(\OO_{\Delta_X},\OO_{\Delta_X}) \longrightarrow
\Ext^*(\OO_{\Delta_S},\OO_{\Delta_S})=:\HH^*(S).
\end{equation}
So although Hochschild cohomology is not in general functorial, we see that it \emph{is} under full and faithful functors (as well as equivalences). For $\HH^2$ we interpret this as saying that a deformation of $D(X)$ induces a deformation of the full subcategory $D(S)$.

We denote the map \eqref{HH*} by $\Phi_P^{\HH^*\!}$, and understand it via the variation of Hodge structure it induces. To access this we first go via Hochschild homology.

\subsection{Hochschild homology} Let $S$ and $X$ have dimensions $m$ and $n$ respectively. Hochschild homology is defined by $\HH_*(S):=\Ext^{m-*}
(\Delta_{S*}\omega_S^{-1},\OO_{\Delta_S})$ and \emph{is} functorial. That is, we get a map $\Phi_P^{\HH_*\!}$
\begin{equation} \label{HH_*}\xymatrix@=15pt{
\Ext^{m-*}(\Delta_{S*}\omega_S^{-1},\OO_{\Delta_S}) \ar[r] \ar@{=}[d] &
\Ext^{n-*}(\Delta_{X*}\omega_X^{-1},\OO_{\Delta_X}) \ar@{=}[d] \\
\HH_*(S) \ar[r]^{\Phi_P^{\HH_*\!}} & \HH_*(X),}
\end{equation}
by the following construction \cite{Ca1}. We modify diagram \eqref{diagram} to
\begin{equation} \label{diagram2} \hspace{-15mm}
\xymatrix@R=2pt{
D(S\times S) \ar[r]^{P*}& D(S\times X) \ar[r]^{*R} & D(X\times X) \\
\qquad \OO_{\Delta_S}\quad \ar@{|->}[r]& \qquad P \qquad \ar@{|->}[r]& \quad P*R\qquad \\
\quad \Delta_{S*}\omega_S^{-1}\,\ \ar@{|->}[r]& \ \,P\otimes\omega_S^{-1}\,\ \ar@{|->}[r]& P\!*\!L\otimes\pi_1^*\omega^{-1}_X[m-n].\hspace{-16mm}}
\end{equation}
Here, if $\tau$ denotes the isomorphism $X\times S\to S\times X$ that interchanges the two factors, then
\begin{equation} \label{Radjo}
R:=\tau^*(P^\vee\otimes\omega_S[m])\in D(X\times S)
\end{equation}
is the kernel for the right adjoint of $\Phi_P$, and
\[ L:=\tau^*(P^\vee\otimes\omega_X[n])=R\otimes\omega_S^{-1}\otimes\omega_X[n-m] \]
is the kernel for the left adjoint. From \eqref{diagram2} we get a map
\[ \Ext^{m-*}(\Delta_{S*}\omega_S^{-1},\OO_{\Delta_S}) \to
\Ext^{m-*}(P\!*\!L\otimes\pi_1^*\omega^{-1}_X[m-n],P\!*\!R). \]
Compose with the natural maps of kernels
\begin{equation} \label{unit}\xymatrix{
\OO_{\Delta_X} \ar[r]^(.45)\eta & P*L, & P*R \ar[r]^\epsilon & \OO_{\Delta_X}}
\end{equation}
that induce the unit and counit of the adjunctions. This takes us to the
group $\Ext^{n-*}(\Delta_{X*}\omega_X^{-1},\OO_{\Delta_X})=\HH_*(X)$, thus defining the map \eqref{HH_*}. \medskip

For any variety $Y$ there is an action of $\HH^*(Y)$ on $\HH_*(Y)$ given by composition:
\begin{equation} \label{HHaction}
\xymatrix@=20pt{
\!\Ext^i(\Delta_{Y*}\omega_Y^{-1},\OO_{\Delta_Y})\otimes
\Ext^j(\OO_{\Delta_Y},\OO_{\Delta_Y}) \ar[r] &
\Ext^{i+j}(\Delta_{Y*}\omega_Y^{-1},\OO_{\Delta_Y}).\!\!\!}
\end{equation}
This action has a certain compatibility with the maps $\Phi_P^{\HH^*\!}$ \eqref{HH*} and $\Phi_P^{\HH_*\!}$ \eqref{HH_*}, as described in the next result.

\begin{prop}
Fix $P\in D(S\times X)$ and $e_X\in\HH^j(X)$. The following diagram commutes:
\begin{equation} \label{dg3}
\xymatrix@C=35pt{
\HH_i(S) \ar[r]^{\Phi_P^{\HH_*\!}}\ar[d]_{\Phi_P^{\HH^*\!}(e_X)} & \HH_i(X)\ar[d]^{e_X} \\
\HH_{i-j}(S) \ar[r]^{\Phi_P^{\HH_*\!}} & \HH_{i-j}(X).\!}
\end{equation}
\end{prop}

\begin{proof}
Fix any $f_S\in\HH_i(S)$ mapping via \eqref{HH_*} to $f_X\in\HH_i(X)$. Consider $f_S$ as a map from the bottom left hand object $\Delta_{S*}\omega_S^{-1}$ in the diagram \eqref{diagram2} to the central left hand object $\OO_{\Delta_S}[m-i]$.
Compose with $\Phi_P^{\HH^*\!}(e_X)\colon\OO_{\Delta_S}\to\OO_{\Delta_S}$ and follow the composition from left to right through \eqref{diagram2}. We get
\[ \xymatrix@=20pt{
\quad\OO_{\Delta_S}\quad & P & P*R \ar[r]^\epsilon & \OO_{\Delta_X} \\
\quad\OO_{\Delta_S}\quad \ar[u]^{\Phi_P^{\HH^*\!}(e_X)} \ar@{|->}^(.45){P*}[r]& \qquad\ P\ \qquad \ar[u]^{e_X*1_P} \ar@{|->}^{*R}<112pt,-37pt>;<142pt,-37pt>& P*R \ar[u]^{e_X*1_P*1_R}\ar[r]^\epsilon & \OO_{\Delta_X} \ar[u]^{e_X} \\
\quad\Delta_{S*}\omega_S^{-1}\quad \ar[u]^{f_S}& P\otimes\omega_S^{-1} \ar[u]& P\!*\!L\otimes\pi_1^*\omega^{-1}_X[n-m] \ar[u] \\
&&\Delta_{X*}\omega_X^{-1}[n-m]. \ar[u]^(.45)\eta} \]
Applying $P*$ to $\Phi_P^{\HH^*\!}(e_X)$ (top left) gives $e_X*1_P$ (top centre) by the definition of $\Phi_P^{\HH^*\!}(e_X)$ via diagram \eqref{diagram}. This explains the labelling of the arrows. We are required to prove that the full composition up-up-up-right on the right hand side is $e_X\circ f_X$. It is the same as up-up-right-up because the top right square clearly commutes. But up-up-right is $f_X$ by definition.
\end{proof}

\subsection{Hodge cohomology}
We can recast the results of the last two sections in terms of Hodge theory via Kontsevich's modification of the Hochschild-Kostant-Rosenberg isomorphisms
\begin{equation} \label{HKR}
\HH^*(Y)\ \cong\ \bigoplus_{i+j=*}H^i(\Lambda^jT_Y),\qquad
\HH_*(Y)\ \cong\ \bigoplus_{j-i=*}H^i(\Omega^j_Y).
\end{equation}
These are defined by post-composing the standard HKR isomorphism (given by the exponential of the universal Atiyah class  \cite{BuFl, Ca2}) with $\td^{-1/2}\hook(\ \cdot\ )$ acting on $\bigoplus H^i(\Lambda^jT_Y)$ and $\td^{1/2}\wedge(\ \cdot\ )$ acting on $\bigoplus H^i(\Omega^j_Y)$. Whenever we write something like $H^i(T_Y)\subset\HH^{i+1}(Y)$ we mean via this modified HKR isomorphism.

For us, the modified HKR isomorphism has two main advantages over its untwisted cousin. Firstly \cite{CV1, CV2}, it intertwines the action \eqref{HHaction} of $\HH^*(Y)$ on $\HH_*(Y)$ with the interior multiplication of $H^*(\Lambda^*T_Y)$ on $H^*(\Omega^*_Y)$.
Secondly \cite[Thm.~1.2]{MS1}, it intertwines the map $\Phi_P^{\HH_*\!}$ \eqref{HH_*} on Hochschild homology with the usual map $\Phi_P^H$ \eqref{HP} on cohomology. \medskip

From now on we restrict attention to the case of interest: $S$ a K3 surface, $X$ a cubic fourfold, and $P\in D(S\times X)$ the kernel of a full and faithful embedding $D(S)\to D(X)$. The holomorphic 2-form (unique up to scale)
\[ \sigma_S\in H^0(\Omega^2_S)\cong\HH_2(S)\cong\C \]
generates $\HH_2(S)$. Similarly let
\[ \sigma_X\in H^1(\Omega^3_X)\cong\HH_2(X)\cong\C \]
be its image under $\Phi_P^H$ (or equivalently $\Phi_P^{\HH_*\!}$) -- also a generator since $\Phi_P^H$ is injective. Using these we can see when a first order commutative deformation of $D(X)$ induces a commutative deformation of $D(S)$ entirely in terms of Hodge theory.

\begin{prop} \label{HHH} Let $\kappa^{}_X\in H^1(T_X)\subset\HH^2(X)$ be a first order deformation of the cubic fourfold $X$, defining a cohomology class
\[ \tilde\kappa^{}_X:=\kappa^{}_X\hook\sigma_X\in H^{2,2}(X)\subset\HH_0(X). \]
Suppose that $\tilde\kappa^{}_X=\Phi_P^H(\tilde\kappa^{}_S)$ for some $(1,1)$ cohomology class $\tilde\kappa^{}_S\in H^{1,1}(S)\subset\HH_0(S)$. Writing
$\tilde\kappa^{}_S=\kappa^{}_S\hook\sigma_S$ for $\kappa^{}_S\in H^1(T_S)$, we have
\[
\Phi_P^{\HH^*\!}(\kappa^{}_X)=\kappa^{}_S.
\]
\end{prop}

\begin{proof}
We follow $\sigma_S\in H^{2,0}(S)\subset\HH_2(S)$ clockwise around the commutative diagram \eqref{dg3} from the top left hand corner:
\[ \xymatrix@C=35pt{
\HH_2(S) \ar[r]^{\Phi_P^{\HH_*\!}}\ar[d]_{\Phi_P^{\HH^*\!}(\kappa^{}_X)} &
\HH_2(X)\ar[d]^{\kappa^{}_X} \\
\HH_0(S) \ar[r]^{\Phi_P^{\HH_*\!}} & \HH_0(X).\!} \]
It maps to $\sigma_X$ in the top right, then $\tilde\kappa^{}_X$ in the bottom right. By assumption this is
$\Phi_P^{\HH_*\!}(\tilde\kappa^{}_S)$, and $\Phi_P^{\HH_*\!}$ is an injection, so on the left hand side we find that $\Phi_P^{\HH^*\!}(\kappa^{}_X)\lrcorner\,\sigma_S
=\kappa^{}_S\lrcorner\,\sigma_S$. 

So the result follows from the fact that contraction with $\sigma_S$ gives an isomorphism $\HH^*(S)\to\HH_{2-*}(S)$ (corresponding via HKR to the isomorphism $H^i(\Lambda^jT_S)\to H^i(\Omega^{2-j}_S)$ induced by $\hook\sigma_S\colon \Lambda^j T_S\xrightarrow{\sim}\Omega^{2-j}_S$). \end{proof}

\section{Deformations} \label{2}

\subsection{First order} \label{2.1}
We now have a Hodge-theoretic criterion for a commutative deformation $\kappa^{}_X\in H^1(T_X)$ to define via $\Phi_P^{\HH^*\!}$ a commutative deformation $\kappa^{}_S\in H^1(T_S)$ of $S$. In this section we will show that to first order, if we deform $X$ and $S$ by these Kodaira-Spencer classes $\kappa^{}_X$ and $\kappa^{}_S$ respectively, the fully faithful Fourier-Mukai kernel $P\in D(S\times X)$ deforms with them. Let $A_n$ denote $\Spec\C[t]/(t^{n+1})$.

\begin{thm} \label{A1}
Suppose that $\kappa^{}_X\in H^1(T_X)\subset\HH^2(X)$ maps via $\Phi_P^{\HH^*\!}$ to $\kappa^{}_S\in H^1(T_S)\subset\HH^2(S)$. Let $X_1\to A_1$ and $S_1\to A_1$ be the corresponding flat deformations of $X$ and $S$. Then there is a perfect complex
\[ P_1\in D(S_1\times_{A_1}X_1) \]
whose derived restriction to $S\times X$ is $P$.
\end{thm}

By \cite[Thm.~3.1]{HMS}, \cite[Cor.~3.4]{HT} we only have to prove the vanishing of the obstruction class
\begin{equation} \label{KSAt}
(\kappa^{}_S,\kappa^{}_X)\circ\At(P)\ \in\,\Ext^2_{S\times X}(P,P).
\end{equation}
Here
\begin{equation} \label{At1}
\At(P) \colon P \longrightarrow P\otimes\Omega_{S\times X}[1]
\end{equation}
is the Atiyah class of $P$. The key to the proof is the interaction of various Atiyah classes with Fourier-Mukai transforms; this is what we turn to in the next section.
We think of the Kodaira-Spencer class $(\kappa^{}_S,\kappa^{}_X)\in H^1(T_{S\times X})$ as a morphism
\[ \Omega_{S\times X} \xrightarrow{(\kappa^{}_S,\kappa^{}_X)} \OO_{S\times X}[1] \]
which we compose with the Atiyah class \eqref{At1} to give the obstruction morphism $P\to P[2]$.

\subsection*{Atiyah classes}
Fix a complex $F\in D(Y)$ on a smooth projective variety $Y$. Recall the first jet space $J^1(F)$ of $F$ is defined by
\begin{equation} \label{jett}
J^1(F):=\pi_{2*}\big(\pi_1^*F\otimes\OO_{2\Delta_Y}\big).
\end{equation}
Here $\pi_i$ is projection onto the $i$th factor of $Y\times Y$, and $2\Delta_Y\subset
Y\times Y$ is the subscheme whose ideal sheaf is the square $\I_{\!\Delta_Y}^2$ of the ideal sheaf of the diagonal. In other words $J^1(F)$ is the image of $F$ under the Fourier-Mukai transform $D(Y)\to D(Y)$ with kernel $\OO_{2\Delta_Y}$. The obvious exact sequence of kernels\footnote{Here and below we identify the conormal bundle to $\Delta_Y$ with $\Omega_Y$ via the \emph{first} projection $\pi_1$. This (arbitrary) choice is the origin of the signs in Corollary \ref{AtAt}, which are flipped if we instead use $\pi_2$. \label{signn}}
\begin{equation} \label{2diag}
0 \longrightarrow \Delta_{Y*}\Omega_Y \longrightarrow \OO_{2\Delta_Y} \longrightarrow \OO_{\Delta_Y} \longrightarrow 0,
\end{equation}
applied to $F$ gives the standard exact triangle
\begin{equation} \label{At}
F\otimes\Omega_Y \longrightarrow J^1(F) \longrightarrow F.
\end{equation}
The Atiyah class $\At(F)\in\Ext^1(F,F\otimes\Omega_Y)$ of $F$ is defined to be the connecting homomorphism $\At(F)\colon F\to F\otimes\Omega_Y[1]$ (or extension class) of \eqref{At}. For this reason, the extension class of \eqref{2diag} is called the \emph{universal Atiyah class}
\begin{equation} \label{UAt}
\mathbf{At}_Y\in\Ext^1_{Y\times Y}(\OO_{\Delta_Y},\Delta_{Y*}\Omega_Y).
\end{equation}

\subsection*{Partial Atiyah classes}
When $Y=A\times B$, the Atiyah class splits,
\begin{eqnarray}
\Ext^1(F,F\otimes\Omega_{A\times B}) &\cong&
\Ext^1(F,F\otimes\Omega_A)\oplus\Ext^1(F,F\otimes\Omega_B), \nonumber \\
\At(F)\qquad &=& \qquad\qquad (\At_A(F)\,,\,\At_B(F)). \label{partial}
\end{eqnarray}
We are suppressing some obvious pullback maps for appearance's sake. This defines the \emph{partial Atiyah class} $\At_A(F)$. (Equivalently it is the relative Atiyah class $\At_{A\times B/B}(F)$.) We describe it by noting that the exact sequence \eqref{2diag} has a natural quotient:
\begin{equation} \label{ABA}
\xymatrix@=15pt{
0 \ar[r] & (\Delta_{A\times B})_*\Omega_{A\times B} \ar[r]\ar[d] & \OO_{2\Delta_{A\times B}} \ar[r]\ar[d] & \OO_{\Delta_{A\times B}} \ar[r]\ar@{=}[d] & 0 \\
0 \ar[r] & (\Delta_{A\times B})_*\Omega_A \ar[r] & \OO_{2\Delta_A\times\Delta_B} \ar[r] & \OO_{\Delta_{A\times B}} \ar[r] & 0,\!}
\end{equation}
where $2\Delta_A\times\Delta_B$ is defined by the ideal sheaf $\I_{\!\Delta_A}^2\boxtimes
\I_{\!\Delta_B}$. The lower row is the extension defined by the extension class of the top row, projected via $\Omega_{A\times B}\to\Omega_A$. Therefore it is what defines $\At_A(F)$. It is the sequence 
\begin{equation} \label{messB}
0 \longrightarrow \Delta_{A*}\Omega_A\boxtimes\OO_{\Delta_B} \longrightarrow \OO_{2\Delta_A}\boxtimes\OO_{\Delta_B} \longrightarrow \OO_{\Delta_A}\boxtimes\OO_{\Delta_B} \longrightarrow 0,
\end{equation}
i.e.\ ${}\boxtimes\OO_{\Delta_B}$ applied to the sequence \eqref{2diag} for $Y=A$. 
Now apply this to $F\in D(A\times B)$. We pull up to $A\times B\times A\times B$ from the first two factors, tensor with \eqref{messB}, and push down to the second two factors. This push down factors through the push down to $A\times B\times A$, which turns \eqref{messB} into
\begin{equation} \label{messB2}
0 \longrightarrow \Delta_{A*}\Omega_A\boxtimes\OO_B \longrightarrow \OO_{2\Delta_A}\boxtimes\OO_B \longrightarrow \OO_{\Delta_A}\boxtimes\OO_B \longrightarrow 0.
\end{equation}
So equivalently we can pull $F$ back to $A\times B\times A$ from the first two factors, tensor with \eqref{messB2} and push down to the last two factors. The upshot is that if we define the partial jet space by
\[ J^1_A(F):=\pi_{23*}\big(\pi_{12}^*F\otimes\OO_{2\Delta_A\times B}) \]
(where the projections are from $A\times B\times A$ to its factors), then $\At_A(F)$ is the connecting homomorphism of the exact triangle
\begin{equation} \label{AtA}
F\otimes\Omega_A \longrightarrow J^1_A(F) \longrightarrow F
\end{equation}
induced from \eqref{messB2}.

\begin{lem} \label{Pstar}
Consider \eqref{2diag} for $Y=S$ to be an exact triangle of Fourier-Mukai kernels
in $D(S\times S)$:
\begin{equation} \label{again}
\Delta_{S*}\Omega_S \longrightarrow \OO_{2\Delta_S} \longrightarrow \OO_{\Delta_S}.
\end{equation}
Apply the Fourier-Mukai functor $P*$ of diagram \eqref{diagram}. Then the resulting exact triangle of kernels in $D(S\times X)$ is the triangle \eqref{AtA} defining $\At_S(P)$:
\begin{equation} \label{AtSP}
P\otimes\Omega_S \longrightarrow J^1_S(P) \longrightarrow P.
\end{equation}
\end{lem}

\begin{proof}
Applying $P*$ means we have to pull \eqref{again} back to $S\times S\times X$,
\begin{equation} \label{sx}
\Delta_{S*}\Omega_S\boxtimes\OO_X \longrightarrow \OO_{2\Delta_S}\boxtimes\OO_X \longrightarrow \OO_{\Delta_S}\boxtimes\OO_X,
\end{equation}
then tensor with $\pi_{23}^*P$ and push down by $\pi_{13*}$.

But \eqref{sx} is precisely the sequence \eqref{messB2} for $A=S,\,B=X$ (with the second and third factors swapped). So by \eqref{AtA} we indeed get the triangle \eqref{AtSP} with extension class $\At_S(P)$.
\end{proof}

Symmetrically we similarly find the following.

\begin{lem} \label{starP}
Consider \eqref{2diag} for $Y=X$ to be an exact triangle of Fourier-Mukai kernels
in $D(X\times X)$:
\[ \Delta_{X*}\Omega_X \longrightarrow \OO_{2\Delta_X} \longrightarrow \OO_{\Delta_X}. \]
Apply the Fourier-Mukai functor $*P$ of diagram \eqref{diagram}. Then the resulting exact triangle of kernels in $D(S\times X)$ is the triangle \eqref{AtA} defining $\At_X(P)$:
\begin{flalign*}
\phantom{\qed} & & P\otimes\Omega_X \longrightarrow J^1_X(P) \longrightarrow P. & &\qed
\end{flalign*}
\end{lem}

\subsection*{The universal Atiyah class}
The Atiyah class $\At(\OO_{\Delta_S})$ of $\OO_{\Delta_S}$ on $S\times S$ lies in $\Ext^1_{S\times S}(\OO_{\Delta_S},\Delta_{S*}(\Omega_{S\times S}|_{\Delta_S}))$, so cannot quite be the universal Atiyah class $\textbf{At}_S\in\Ext^1_{S\times S}(\OO_{\Delta_S},\Delta_{S*}\Omega_S)$ of \eqref{UAt}. To see their close relationship \eqref{tues} we need the following standard fact.

\begin{lem} Let $Z\subset Y$ be a subscheme of a smooth variety, giving the standard short exact sequence
\begin{equation} \label{2times}
0\longrightarrow\I_Z/\I_Z^2\longrightarrow\OO_{2Z}\longrightarrow\OO_Z
\longrightarrow0.
\end{equation}
The image of its extension class in $\Ext^1(\OO_Z,\I_Z/\I_Z^2)$ under the canonical map $\I_Z/\I_Z^2\to\Omega_Y|_Z$ is the Atiyah class
\[ \At(\OO_Z)\in\Ext^1(\OO_Z,\OO_Z\otimes\Omega_Y). \]
\end{lem}

\begin{proof}
We apply $\pi_{2*}(\pi_1^*\OO_Z\otimes\ \cdot\ )$ to the exact sequence
\[ 0\longrightarrow\Omega_Y\otimes\OO_{\Delta_Y}\longrightarrow\OO_{2\Delta_Y}
\longrightarrow\OO_{\Delta_Y}\longrightarrow0 \]
on $Y\times Y$. Since these sheaves are $\pi_1$-flat, we can take the underived tensor product of structure sheaves, giving the structure sheaves of the intersections. The result is $\pi_{2*}$ of the exact sequence
\[ 0\longrightarrow\Omega_Y\otimes\OO_{\Delta_Z}\longrightarrow\OO_{2(\Delta_Z\subset Z\times Y)} \longrightarrow\OO_{\Delta_Z}\longrightarrow0. \]
Here $2(\Delta_Z\subset Z\times Y)$ denotes the doubling of $\Delta_Z$ inside $Z\times Y$. This scheme maps onto $2Z\subset Y$ by the projection $\pi_2$. Therefore using $\pi_2^*$ to pull back functions gives a map of sheaves from $\OO_{2Z}$ to $\OO_{2(\Delta_Z\subset Z\times Y)}$. Of course this is \emph{not} a map of $\OO_{Y_1\times Y_2}$-modules, but it \emph{is} a map of $\OO_{Y_2}$-modules. Therefore after we apply $\pi_{2*}$ we get a map of $\OO_Y$-modules
\[ \xymatrix@=16pt{
0 \ar[r] & \Omega_Y\otimes\OO_Z \ar[r] & \pi_{2*}\OO_{2(\Delta_Z\subset Z\times Y)} \ar[r] & \OO_Z \ar[r] & 0 \\
0 \ar[r] & \I_Z/\I_Z^2 \ar[r] & \OO_{2Z} \ar[r]\ar[u] & \OO_Z \ar[r] & 0.\!}
\]
It is easy to see that this completes in the obvious way with vertical arrows on the left and right.
Since the upper central term is $J^1(\OO_Z)$ \eqref{jett} and the upper row is the sequence \eqref{At} defining $\At(\OO_Z)$, this gives the result claimed.
\end{proof}

We apply this to $Z=\Delta_S$ inside $Y=S\times S$, for which the canonical map $\I_Z/\I_Z^2\to\Omega_Y|_Z$ splits:
\begin{equation} \label{splat}
\Omega_{S\times S}|_{\Delta_S}\ \cong N^*_{\Delta_S}\oplus\Omega_{\Delta_S}.
\end{equation}
Here $N^*_{\Delta_S}$ denotes the conormal bundle. Furthermore, the sequence \eqref{2times} becomes \eqref{2diag}, with extension class the universal Atiyah class $\mathbf{At}_S$ of \eqref{UAt}. Therefore $\At(\OO_{\Delta_S})$ splits too, with respect to \eqref{splat}, as
\begin{equation} \label{tues}
\At(\OO_{\Delta_S})=\mathbf{At}_S\oplus0\ \in\,\Ext^1\big(\OO_{\Delta_S},\OO_{\Delta_S}
\otimes (N_S^*\oplus\Omega_{\Delta_S})\big).
\end{equation}
In particular we have proved the following.

\begin{cor} \label{AtAt}
Denote the two partial Atiyah classes \eqref{partial} of $\OO_{\Delta_S}$ by \newline $\At_1(\OO_{\Delta_S})$ and $\At_2(\OO_{\Delta_S})$. Then (cf. footnote \ref{signn}),
\[
\At_1(\OO_{\Delta_S})=\mathbf{At}_S=-\At_2(\OO_{\Delta_S})\ \in\
\Ext^1(\OO_{\Delta_S},\Omega_{\Delta_S}).
\]
\end{cor}
\smallskip

\begin{proof}[Proof of Theorem \ref{A1}.]
Following Toda \cite{Toda}, the idea of the proof is to see the HKR isomorphism as identifying: (1) first order deformations of $D(S)$, with (2) the corresponding obstructions to deforming the identity Fourier-Mukai kernel $\OO_{\Delta_S}\in D(S\times S)$ as the first $S$ (or $D(S)$) factor deforms and the second remains fixed. Applying that same philosophy to $P\in D(S\times X)$ as well quickly gives the result.

The Kodaira-Spencer class $\kappa^{}_S\in H^1(T_S)$ induces a corresponding first order deformation $\pi_1^*\kappa^{}_S\in H^1(T_{S\times S})$ of $S\times S$ (deforming the first factor and fixing the second).
The obstruction to deforming $\OO_{\Delta_S}\in D(S\times S)$ with this deformation is \cite{HMS,HT}
\[ \pi_1^*\kappa^{}_S\circ\At(\OO_{\Delta_S})\in\Ext^2(\OO_{\Delta_S},\OO_{\Delta_S})=\HH^2(S). \]
By \eqref{tues} and Corollary \ref{AtAt} this equals
\begin{equation} \label{kaps}
\pi_1^*\kappa^{}_S\circ\At_1(\OO_{\Delta_S})=\pi_1^*\kappa^{}_S\circ\textbf{At}_S, \end{equation}
which is precisely the image of $\kappa^{}_S$ under the HKR
isomorphism \eqref{HKR}. (This follows from two simplifications. Firstly the $\td^{-1/2}$ twisting acts as the identity on $H^1(T_Y)$ because the map $H^1(T_Y)\to H^2(\OO_Y)$ given by contraction with $\td_1^{-1/2}=-c_1(Y)/4$ calculates the derivative of the $(0,2)$-Hodge piece of $c_1(\omega_Y)/4$ as $Y$ is deformed by an $H^1(T_Y)$ class, which is of course zero. Secondly,
in this low degree,
the map \,$\circ\exp(\textbf{At}_S)$ of \cite{BuFl,Ca2} reduces to
$\circ\textbf{At}_S\colon H^1(T_S)\hookrightarrow\Ext^2(\OO_{\Delta_S},
\OO_{\Delta_S})=\HH^2(S)$.)

Similarly, the image of $\kappa^{}_X\in H^1(T_X)$ in $\Ext^2(\OO_{\Delta_X},
\OO_{\Delta_X})=\HH^2(X)$ is given by
\[ \pi_1^*\kappa^{}_X\circ\At_1(\OO_{\Delta_X})=-\pi_2^*\kappa^{}_X\circ\At_2(\OO_{\Delta_X}) \]
by Corollary \ref{AtAt}. To describe its image under $\Phi_P^{\HH^*\!}$ \eqref{HH*} we have to chase around the diagram \eqref{diagram}. Under the Fourier-Mukai composition $*P$ of that diagram it maps to 
\begin{equation} \label{Xobs}
-\pi_2^*\kappa^{}_X\circ\At_X(P)\in\Ext^2_{S\times X}(P,P)
\end{equation}
by Lemma \ref{starP}.

Similarly, under the Fourier-Mukai composition $P*$ of diagram \eqref{diagram}, $\kappa^{}_S$ (or rather its image \eqref{kaps} in $\HH^2(S)$) is taken to
\begin{equation} \label{Sobs}
\pi_1^*\kappa^{}_S\circ\At_S(P)\in\Ext_{S\times X}^2(P,P)
\end{equation}
by Lemma \ref{Pstar}.

Since $\Phi_P^{\HH^*\!}(\kappa^{}_X)=\kappa^{}_S$, the classes (\ref{Sobs}, \ref{Xobs}) are equal. In particular,
\[ (\kappa^{}_S,\kappa^{}_X)\circ\At(P)=\pi_1^*\kappa^{}_S\circ\At_S(P)+\pi_2^*\kappa^{}_X\circ\At_X(P)
=0 \]
in $\Ext^2(P,P)$. Thus the obstruction class \eqref{KSAt} vanishes as required. \end{proof}

\subsection{All orders by \texorpdfstring{$T^1$}{T1}-lifting} \label{2.2}
We now extend Theorem \eqref{A1} to all orders using $T^1$-lifting methods \cite{Ra, Ka}. The philosophy is that, to extend an equivalence from a family over $A_n$ to one over $A_{n+1}$, it is enough to understand deformations ``sideways" from $A_n$ to $A_n\times A_1$. For this latter problem we can invoke Theorem \ref{A1} on deforming over $A_1$.

We use the standard Artinian spaces
\begin{align*}
A_n &:= \Spec\C[t]/(t^{n+1}), \\
B_n &:= A_n\times A_1=\Spec\C[x,y]/(x^{n+1},y^2), \\
C_n &:= \Spec\C[x,y]/(x^{n+1},xy^n,y^2).
\end{align*}
 The key to $T^1$-lifting is the following maps between them:
\begin{equation} \label{maps}
\xymatrix{
B_{n-1}\,\ar@{^(->}[r]\ar@{->>}[dr]_{p_{n-1}} &
C_n\,\ar@{^(->}[r]\ar@{->>}[d]^{q_n} & B_n \ar@{->>}[d]^{p_n}& x+y \\
& A_n\,\ar@{^(->}[r]^(.45){\iota_n} & A_{n+1},\! & t.\!\! \ar@{|->}[u]}
\end{equation}
The surjections $p_{n-1},q_n$ and $p_n$ are defined by their action on $t$, which they pull back to $x+y$. So they act ``diagonally", mixing up the $A_n$ and $A_1$ axes; we picture them as follows (for $n=1$).

\begin{center}
\begin{picture}(0,0)
\includegraphics{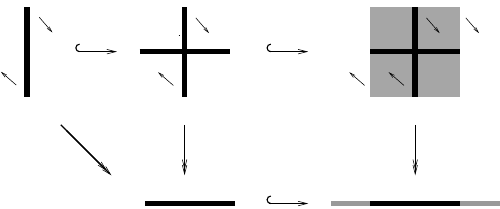}
\end{picture}%
\setlength{\unitlength}{2368sp}%
\begin{picture}(6662,2755)(3390,-5301)
\fontsize{8}{9.6pt} \selectfont
\put(4201,-4561){$p^{}_0$}
\put(5926,-4561){$q^{}_1$}
\put(9001,-4561){$p^{}_1$}
\put(7141,-5086){$\iota^{}_1$}
\end{picture}\end{center}
\medskip

Since $p_n^*\big(\frac{t^{n+1}}{n+1}\big)=x^ny$, the ideal $(t^{n+1})$ of $A_n\subset A_{n+1}$ is isomorphic, via the pull back $p_n^*$, to the
ideal $(x^ny)$ of $C_n\subset B_n$. (They are the ideals of the grey shaded areas in the figure.) Therefore we will find that \emph{extending a deformation from $C_n$ to $B_n$ becomes the same problem as extending from $A_n$ to $A_{n+1}$}. In fact we have the diagram of exact sequences
\begin{equation} \label{T1dg}
\xymatrix@=16pt{(t^{n+1}) \ar[r]\ar@{=}[d] & \OO_{A_{n+1}} \ar[r]\ar[d]^{p_n^*} & \OO_{A_n} \ar[d]^{q_n^*}\\
p_{n*}(x^ny) \ar[r] & p_{n*}\OO_{B_n} \ar[r]\ar[d] & q_{n*}\OO_{C_n} \ar[d] \\ & Q_n \ar@{=}[r] & Q_n}
\end{equation}
on $A_{n+1}$. Here the maps $p_n^*,\,q_n^*$ are $\OO_{A_{n+1}}\to
p_{n*}p_n^*\OO_{A_{n+1}}$ and $\OO_{A_n}\to q_{n*}q_n^*\OO_{A_n}$ respectively. We denote the cokernel of the latter by $Q_n$, and suppress pushforwards by obvious inclusions. We consider \eqref{T1dg} as showing how to recover the flat deformation $\OO_{A_{n+1}}$ over $A_{n+1}$ of $\OO_{A_n}$ over $A_n$ from the flat deformation $\OO_{B_n}$ of $\OO_{C_n}$: i.e.\ push down and take the kernel of the map to $Q_n$. \medskip

The case of deformations of complexes of sheaves is not much harder. Suppose we are given a smooth projective family $\A_{n+1}\to A_{n+1}$. 
Base-changing by $\iota_n$ gives a family $\A_n\to A_n$. Similarly base-changing by $p_n$ and $q_n$ gives $\B_n\to B_n$ and $\cC_n\to C_n$ respectively. We use the same symbols $p_n,q_n,\iota_n$ to denote the induced maps between them.

We start with a perfect complex $P_n\in
D(\A_n)$, whose (derived) restriction to $A_{n-1}$ we denote $P_{n-1}:=\iota_{n-1}^*P_n$.
To extend $P_n$ to some $P_{n+1}$ on $\A_{n+1}$, we consider sideways first order deformations of $P_n$ over $A_n\times A_1=B_n$, i.e.\ perfect complexes on $\B_n$ whose base-change by $A_n\times\{0\}\hookrightarrow A_n\times A_1$ is $P_n$.

\begin{prop} \label{T1prop}
\emph{($T^1$-lifting for complexes of sheaves.)}\footnote{The conditions of this Proposition can be stated more attractively when the family $\A_n$ is a product $Y\times A_n$, i.e.\ when we are not deforming the underlying variety $Y$. We ask for a deformation class $e_{n+1}\in\Ext^1_{\A_n}(P_n,P_n)$ whose restriction to $\A_{n-1}$ gives the extension class
$e_n\in\Ext^1_{\A_{n-1}}(P_{n-1},P_{n-1})$ of $p_{n-1}^*(P_n)$. (The complex $p_{n-1}^*(P_n)$ on $\B_{n-1}$ is a first order deformation of $P_{n-1}$ since that is what it restricts to on $\A_{n-1}$.)}
Suppose there exists a first order deformation
\[ \widetilde P_{n+1}\in D(\B_n) \]
of $P_n\in D(\A_n)$ whose restriction to $\B_{n-1}$ is
\begin{equation} \label{pn-1}
p_{n-1}^*(P_n)\in D\big(\B_{n-1}\big).
\end{equation}
Then there exists $P_{n+1}\in D(\A_{n+1})$ whose restriction to $\A_n$ is $P_n$ (and such that $p_n^*P_{n+1}=\widetilde P_{n+1}$).
\end{prop}

\begin{proof}
By assumption, $\widetilde P_{n+1}$ restricts to $P_n$ on $\A_n$ and to $p_{n-1}^*P_n$ on $\B_{n-1}$, both of which restrict to $P_{n-1}$ on $\A_{n-1}$. 
But $\cC_n=\B_{n-1}\cup_{\A_{n-1}}\A_n$, so the restriction of $\widetilde P_{n+1}$ to $\cC_n$ is $q_n^*P_n$. Suppressing pushforwards by some obvious inclusions, we get the exact triangle
\[ P_0\otimes(x^ny) \longrightarrow \widetilde P_{n+1} \longrightarrow q_n^*P_n \]
on $\B_n$. Now mimic diagram \eqref{T1dg}. Pushing down by $p_{n*}$ gives the middle row of the following diagram of exact triangles on $\A_{n+1}$:
\begin{equation} \label{video}
\xymatrix@=16pt{P_0\otimes(t^{n+1}) \ar[r]\ar@{=}[d] & P_{n+1} \ar[r]\ar[d] & P_n \ar[d]^{q_n^*} \\
P_0\otimes p_{n*}(x^ny) \ar[r] & p_{n*}\widetilde P_{n+1} \ar[r]\ar[d] &
q_{n*}q_n^*P_n \ar[d] \\ & \ Q_n \ar@{=}[r] & \,Q_n.\!\!}
\end{equation}
The right hand column defines $Q_n:=\mathrm{Cone}\,(P_n\to q_{n*}q_n^*P_n)$, then we define $P_{n+1}:=\mathrm{Cone}\,(p_{n*}\widetilde P_{n+1}\to Q_n)[-1]$ by the middle column. \medskip

The arrows $P_{n+1}\to\iota_{n*}P_n$ and $P_{n+1}\to p_{n*}\widetilde P_{n+1}$ in \eqref{video} induce, by adjunction,
\begin{equation} \label{isos}
\xymatrix{\iota_n^*P_{n+1} \ar[r] & P_n & \text{and} &
p_n^*P_{n+1} \ar[r] & \widetilde P_{n+1}.}
\end{equation}
We claim that these are quasi-isomorphisms, which we can check locally as follows.
Since $\widetilde P_{n+1}$ is perfect we may take it to be a finite complex of locally free sheaves whose restriction to $\A_n$ is a finite complex of locally frees representing $P_n$. Locally free sheaves are adapted to all of the pullback and (finite!) pushforward functors in \eqref{video}: they can be applied term by term to the sheaves in the complex to produce the derived functors. Working locally then, we need only check that the maps \eqref{isos} are isomorphisms for free sheaves, which follows from \eqref{T1dg}.

Finally we claim that $P_{n+1}$ is perfect. By the first isomorphism of \eqref{isos}, the (derived) restriction of $P_{n+1}$ to any point of $\A_0$ has cohomology in only finitely many degrees (because $P_n$ is perfect). Therefore, by the Nakayama lemma, an infinite locally free resolution of $P_{n+1}$ can be trimmed to give a finite one.
\end{proof}

In our application we take any smooth curve through the point $0\in\cC_d^\lev$ of Proposition \ref{state}, and complete at $0$. Pulling back $\cS,\,\X$ we get smooth families (also denoted $\cS,\,\X$) of K3 surfaces and cubic fourfolds respectively over the formal curve $A_\infty:=\Spec\C[\![t]\!]$. Let $S_n,\,X_n$ be their restrictions to $A_n\subset A_\infty$. Over the central fibres $S_0$ and $X_0$ we have the kernel
\[ P_0\in D(S_0\times X_0) \]
of a \emph{fully faithful} embedding $D(S_0)\hookrightarrow D(X_0)$, inducing an embedding
\[ \xymatrix{\Phi_{P_0}^H\colon H^*(S_0)\ \ar@{^(->}[r] & H^*(X_0).} \]
Using the natural trivialisations \cite[Prop.~3.8]{Blo}
\begin{equation} \label{GaussManin}
H^*_{dR}(\cS/A_\infty)\cong H^*(S_0)\otimes\C[\![t]\!], \qquad
H^*_{dR}(\X/A_\infty)\cong H^*(X_0)\otimes\C[\![t]\!],
\end{equation}
we get the inclusion
\begin{equation} \label{PH}\xymatrix{
H^*_{dR}(\cS/A_\infty)\ \ar@{^(->}[r] & H^*_{dR}(\X/A_\infty).}
\end{equation}
Up to units in $\C[\![t]\!]$ there is a unique fibrewise holomorphic 2-form
on $\cS/A_\infty$ and a unique fibrewise $(3,1)$-form on $\X/A_\infty$:
\[ \sigma^{}_{\!\cS}\in H^2_{dR}(\cS/A_\infty), \qquad \sigma^{}_{\!\X}\in H^4_{dR}(\X/A_\infty). \]
The span of $\sigma^{}_{\!\cS}$ is taken to the span of $\sigma^{}_{\!\X}$ by $\Phi_{P_0}^H$. Rescaling $\sigma^{}_{\!\X}$ if necessary, we are in the following situation.

\begin{thm} \label{lasttt}
Suppose that $\Phi_{P_0}^H\otimes1_{\C[\![t]\!]}$ maps
\begin{equation} \label{assump}
\xymatrix{\sigma^{}_{\!\cS} \ar@{|->}[r] & \sigma^{}_{\!\X}.}
\end{equation}
Then $P_0$ extends uniquely to all orders. That is, there exist unique $P_n\in D(S_n\times_{A_n} X_n)$ satisfying $\iota_{n-1}^*P_n\cong P_{n-1}$ for all $n\ge1$, defining full and faithful embeddings $\Phi_{P_n}\colon D(S_n)\to D(X_n)$.
\end{thm}

\begin{proof}
We will apply Proposition \ref{T1prop} on the spaces
\[ \A_n:=\cS\times_{A_n}\X \longrightarrow A_n. \]
Base-changing by the maps \eqref{maps} gives similar families $\cC_n\to C_n$ and $\B_n\to B_n$. (Here we also allow $n=\infty$.)

By induction we suppose we have a perfect complex $P_n\in
D(\A_n)$ inducing a fully faithful embedding $D(S_n)\hookrightarrow D(X_n)$. This is certainly true in the base case $n=0$. To produce $P_{n+1}$ we now proceed exactly as in Sections \ref{1} and \ref{2.1}, except now relative to the base $A_n$ in place of $A_0=\Spec\C$.

Considering $\B_n\to A_n\times A_1$ to be a first order deformation (in the category of schemes over $A_n$) of $\A_n\to A_n$, it is characterised by its Kodaira-Spencer class
\[ \kappa_n\in H^1(T_{\A_n/A_n}). \]
With respect to the splitting $\A_n=S_n\times_{A_n}X_n$, this is
\[ \kappa_n=(\kappa^{}_{S_n},\kappa^{}_{X_n})\in H^1\big(T_{S_n/A_n}\oplus T_{X_n/A_n}\big). \]
Contracting with $\sigma^{}_{\!\cS}|^{}_{S_n}\oplus\sigma^{}_{\!\X}|^{}_{X_n}$ gives the relative cohomology class
\[ (\tilde\kappa^{}_{S_n},\tilde\kappa^{}_{X_n})\in H^1\big(\Omega^1_{S_n/A_n}\big)\oplus
H^2\big(\Omega^2_{X_n/A_n}\big). \]
Differentiating \eqref{assump} with respect to $t$ we find that $\Phi_{P_0}^H$ maps
\[ \xymatrix{
H^*(S_0)\otimes\C[\![t]\!]\ \ni\ \dot\sigma^{}_{\!\cS} \ar@{|->}[r] & \dot\sigma^{}_{\!\X}\ \in\ H^*(X_0)\otimes\C[\![t]\!].} \]
By Griffiths' classic variation of Hodge structure calculation\footnote{Griffiths shows that if $\sigma_t\in H^{p,q}(Y_t)\ \forall t$ then
$[\dot\sigma_t]^{p-1,q+1}=\kappa^{}_t\hook\sigma_t\in H^{p-1,q+1}(Y_t)$, where $\kappa^{}_t\in H^1(T_{Y_t})$ is the Kodaira-Spencer class, and we differentiate in the trivialisation \eqref{GaussManin} -- i.e.\ with respect to the Gauss-Manin connection. Projecting $\dot\sigma^{}_{\!\cS}$ to its $(1,1)$ component, and $\dot\sigma^{}_{\!\X}$ to its $(2,2)$ component, removes multiples of $\sigma^{}_{\!\cS}$ and $\sigma^{}_{\!\X}$ respectively. Since these are mapped to each other by \eqref{assump}, this means the projections are too.}, this gives
\[ \xymatrix{
H^*(S_0)\otimes\C[\![t]\!]\ \ni\ \tilde\kappa^{}_{S_\infty} \ar@{|->}[r] &
\tilde\kappa^{}_{X_\infty}\ \in\ H^*(X_0)\otimes\C[\![t]\!].} \]
Base-changing back to $A_n$ we find that \eqref{PH} takes $\tilde\kappa^{}_{S_n}$ to $\tilde\kappa^{}_{X_n}$. Thus by Proposition \ref{HHH} over the base $A_n$, the map on Hochschild cohomology induced by $P_n$ takes $\kappa^{}_{X_n}$ to $\kappa^{}_{S_n}$:
\[ \xymatrix@R=5pt{
\HH^2(\X_n/A_n) \ar[r] & \HH^2(\cS_n/A_n) \\
\qquad\kappa^{}_{X_n}\qquad \ar@{|->}[r] & \qquad\kappa^{}_{S_n}.\qquad} \]
Thus by Theorem \ref{A1} applied to $P_n\in D(\A_n)$ over the base $A_n$ we find that there is a perfect complex $\widetilde P_{n+1}\in D(\B_n)$ whose restriction to $\A_n$ is $P_n$. It is unique since
\[ \Ext^1_{\A_n}(P_n,P_n)\ \stackrel{\eqref{SP}}=\
\HH^1(S_n/A_n)\ \stackrel{\eqref{HKR}}=\ H^1(\OO_{S_n/A_n})\oplus H^0(T_{S_n/A_n}) \]
vanishes for $S_n/A_n$ a family of K3 surfaces.
By uniqueness its restriction to $\B_{n-1}$ is therefore $p_{n-1}^*P_n$, so we can apply Proposition \ref{T1prop} to produce $P_{n+1}\in D(\A_{n+1})$.

The Fourier-Mukai composition of the kernel $P_{n+1}$ with the kernel \eqref{Radjo} of its right adjoint gives a perfect complex on $S_{n+1}\times_{A_{n+1}}S_{n+1}$ whose restriction to $S_0\times S_0$ is $\OO_{\Delta_{S_0}}$. But $\OO_{\Delta_{S_0}}$ is rigid: $\Ext^1_{S_0\times S_0}(\OO_{\Delta_{S_0}},\OO_{\Delta_{S_0}})=0$. Thus the composition is $\OO_{\Delta_{S_{n+1}}}$, so $\Phi_{P_{n+1}}$ is full and faithful.
\end{proof}

\subsection{Proof of Theorem \ref{MAIN}} \label{finally}
Recall the kernel $P_0$ given to us by Proposition \ref{state}. By Theorem \ref{lasttt} it deforms to all orders along any smooth curve in the smooth space $\cC_d^\lev$. It follows by standard deformation theory that its deformations are unobstructed, and it deforms to all orders over $\cC_d^\lev$ itself. Therefore by \cite[Prop.~3.6.1]{Lieblich} the kernel $P_0$ in fact deforms over the formal neighbourhood $\widehat Z:=\Spec\widehat\OO_{\cC_d^{\lev\!},0}$ of our point $0\in\cC_d^\lev$ to give a kernel $P_{\widehat Z}$ in $D(\cS\times_{\widehat Z}\X)$. 

Lieblich \cite{Lieblich} shows that the stack of complexes with no negative self-Exts on the fibres of $\cS\times_{\cC_d^\lev}\X\to\cC_d^\lev$ is a locally finitely presented Artin stack $\M/\cC_d^\lev$. 
Since having no negative self-Exts is an open condition satisfied by $P_0$ \eqref{SP}, $P_{\widehat Z}$ defines a section of $\M\times_{\widehat Z}\cC_d^\lev\to
\widehat Z$. Furthermore, $P_0$ is rigid (another open condition), so this section is in fact the formal neighbourhood of $P_0$ in $\M$.

As $\M$ is locally finite presented, it follows that there exists a smooth pointed scheme $(Z,0)$ mapping to $(\M,P_0)$ taking the formal neighbourhood of $0\in Z$ isomorphically to $\widehat Z\subset\M$. Passing to a Zariski open around $0 \in Z$ we can ensure the projection $Z\to\cC_d^\lev$ is an \'etale map (as it is at 0) onto a Zariski open subset.

The universal complex pulls back to a twisted complex $P$ on $\cS\times_{Z}\X$ whose restriction to any fibre $S\times X$ is an untwisted complex. Arguing as in the proof of Theorem \ref{lasttt}, its composition with its right adjoint
gives a kernel whose restriction to $S_0\times S_0$ is $\OO_{\Delta_{S_0}}$. But $\OO_{\Delta_{S_0}}$ is rigid, so its restriction to $S_t\times S_t$ is $\OO_{\Delta_{S_t}}$ for any $t$ in a Zariski open about $0$.

Shrinking $Z$ further if necessary, the resulting full and faithful embeddings $D(S_t)\to D(X_t)$ have image in $\A_{X_t}$ because this is true at $t=0$ and being orthogonal to the pullbacks from $X_t$ of $\OO_{X_t},\OO_{X_t}(1),\OO_{X_t}(2)$ is also an open condition.  Since the Serre functor on both $D(S_t)$ and $\A_{X_t}$ is the shift by $[2]$, and since $D(S_t)$ is non-zero and $\A_{X_t}$ is indecomposable, any fully faithful functor $D(S_t) \to \A_{X_t}$ is necessarily an equivalence \cite[Prop.~1.54]{huybrechts_fm}.

Therefore projecting $Z$ by the finite map $\cC_d^\lev\to\cC_d$ gives the Zariski open subset claimed in Theorem \ref{MAIN}.

\subsection{Outlook} \label{end}
We consider the problem of extending Theorem \ref{MAIN} to the whole of $\cC_d$ to be a purely technical problem, but one for which current techniques seem inadequate. The problem is to produce from a family of equivalences $\Phi_t\colon D(S_t)\to\A_{X_t}$ a limiting equivalence
$\Phi_\infty\colon D(S_\infty)\to\A_{X_\infty}$, where $S_\infty$ is the K3 surface associated to $X_\infty$ by Hassett.

The quickest route would be to use stability conditions of some kind, if they were better understood. Ideally we would see $S_t$ as a moduli space of \emph{stable} objects in $\A_{X_t}\subset D(X_t)$ and pass to the limit $t=\infty$ this way. Alternatively we would like a stability condition on $S_t\times X_t$ in which the Fourier-Mukai kernel $P_t$ of $\Phi_t$ is stable, and take the limit of $P_t$ as a stable object of $D(S_\infty\times X_\infty)$.
\smallskip

An alternative is to switch points of view from considering $S_t$ to be a moduli space of objects in $\A_{X_t}\subset D(X_t)$ and instead consider $X_t$ as parametrising objects on $S_t$. Projecting $\OO_x,\ x\in X_t$, into $\A_{X_t}$ gives an object $P_x$ which sits inside an 8-dimensional holomorphic symplectic moduli space $M_t$ of objects of $\A_{X_t}\subset D(X_t)$, as studied in the recent preprint \cite{lehn}. Using the equivalence $\Phi_t$ we can also consider $M_t$ to be a moduli space of objects on the K3 surface $S_t$; when we do this we emphasise this point of view by calling it $N_t\cong M_t$.  The map \[ f_t\colon X_t\longrightarrow N_t,\qquad x\mapsto P_x, \]
expresses $X_t$ as a complex Lagrangian submanifold of $N_t$. Since K3 surfaces are so well understood, taking the limit moduli space $N_\infty$ of stable objects on $S_\infty$ should be no problem, so the issue becomes how to take the limit of the maps $f_t$.

\emph{A priori} this might involve blowing up $X_\infty$, so instead we might try to take the limit of the isomorphisms $M_t\cong N_t$, at least near $X_t\subset M_t$. Namely, replacing $M_t$ by a Zariski open neighbourhood of $X_t\subset M_t$ gives a Zariski open in the Artin stack of \emph{all} objects (with no negative self-Exts) of $D(X_t)$. This stack \emph{does} behave well in families \cite{Lieblich}, giving a limiting stack of objects in $D(X_\infty)$. These include $P_x,\ x\in X_\infty$, so this stack contains the scheme $X_\infty$, and we can set $M_\infty$ to be a Zariski open neighbourhood of this. The upshot is two families of quasi-projective \emph{holomorphic symplectic} varieties $M_t,\ N_t,\ t\in\C\cup\{\infty\}$, isomorphic away from $t=\infty$. We expect that this gives a birational equivalence between $M_\infty$ and $N_\infty$, and so a derived equivalence between compactly supported derived categories too.  The upshot should be a kernel on $S_\infty\times X_\infty$ which we would expect to define a fully faithful embedding. There is clearly a lot of work to do here. \smallskip

Finally we could instead try to use $F(X_t)$ \cite{BD}, the Fano variety of lines on $X_t$, which is also a moduli space of objects in $\A_{X_t}$ (namely the projection into $\A_{X_t}$ of the structure sheaves of the lines \cite[\S5]{km}). The equivalence $\Phi_t$ therefore makes it isomorphic to a 4-dimensional holomorphic symplectic moduli space of objects in $D(S_t)$; in this guise we call it $\mathcal M_t$. In the limit the isomorphism $F(X_t)\cong \mathcal M_t$ makes $F(X_\infty)$ and $\mathcal M_\infty$ birational and so derived equivalent by work of Kawamata and Namikawa.

Then consider the composition
\begin{equation} \label{markma}
D(S_\infty)\longrightarrow D(\mathcal M_\infty)\longrightarrow D(F(X_\infty))\longrightarrow \A_{X_\infty}\subset D(X_\infty),
\end{equation}
where the first arrow is given by the universal complex, the second is the derived equivalence, and the third uses the universal complex that sees $F(X_\infty)$ as a moduli space of objects in $\A_{X_\infty}$. We expect this composition to be of the form $\Phi_\infty\oplus\Phi_\infty[-2]$, where $\Phi_\infty$ is the equivalence we seek.

One way to try to prove this is as follows. Consider the (adjoint of) the composition of the first two arrows as expressing $F(X_\infty)$ as a moduli space of objects on $S_\infty$. One should check that they are \emph{simple} objects. They form a spanning set, so it should be enough to see that the composition \eqref{markma} acts as a direct sum $\Phi_\infty\oplus\Phi_\infty[-2]$ of an equivalence and its shift on these objects. So we reduce to studying the composition
\[ D(F(X_\infty))\longrightarrow D(S_\infty)\longrightarrow D(F(X_\infty))\longrightarrow \A_{X_\infty}\subset D(X_\infty). \]
Now the composition of the first two arrows is the endofunctor of $D(F(X_\infty))$ studied by Markman and Mehrotra in \cite{mm_in_prep}. Though the definition of this endo\-functor involves seeing $F(X_\infty)$ as a moduli space of objects on $S_\infty$, they expect it to be independent of this description, and to be canonically associated to the holomorphic symplectic manifold $F(X_\infty)$; indeed, this would follow from their conditional result \cite[Thm.~1.11]{mm_conditional}.  Thus we should get the same endofunctor by thinking of $F(X_\infty)$ as a moduli space of objects of $\A_{X_\infty}$ via the (adjoint of) the third arrow above, and the composition becomes
\[ D(F(X_\infty))\longrightarrow\A_{X_\infty}\longrightarrow D(F(X_\infty))\longrightarrow \A_{X_\infty}\subset D(X_\infty). \]
But by \cite[Thm.~4]{nick}, the composition of the second and third arrows in this last sequence is $\operatorname{id} \oplus [-2]$, as required.

%% file: hodge_conj.tex

\section{Algebraic cycles} \label{hodge_conj}

In this final section we prove Theorems \ref{algcycle} and \ref{HODGE} from the Introduction.

\begin{proof}[Proof of Theorem \ref{algcycle}]
We use the families $\X\to\cC_d^\lev,\ \cS\to\cC_d^\lev$ of Proposition \ref{state} and the twisted complex $P$ on $\cS\times_{Z}\X$ of Section \ref{finally}, where $Z\subset\cC_d^\lev$ is Zariski open.

The restriction of $P$ to any fibre $S_t\times X_t$ is the (untwisted) kernel of a Fourier-Mukai equivalence $D(S_t)\to\A_{X_t}$. Consider its Mukai vector to be a rational cycle and take its 3-dimensional part. By adding a large multiple of $(\ell+h)^3$ we may assume it is effective for all $t\in Z$. This does not affect the fact that it induces a Hodge isometry $H^2_\prim(S_t,\Z)(-1)\to\langle h^2, T \rangle^\perp$. 

Therefore the appropriate relative Hilbert scheme of subschemes of the fibres of $\cS\times_{\cC_d^\lev}\X\to\cC_d^\lev$ has non-empty fibres over $Z\subset\cC_d^\lev$. Since it is proper, it has non-empty fibres over all of $\cC_d^\lev$. In particular over a point corresponding to $X$ we find a polarised K3 surface $(S,\ell)$ and a cycle in $S\times X$ inducing the sought-after Hodge isometry
$H^2_\prim(S,\Z)(-1)\to\langle h^2, T \rangle^\perp$.
\end{proof}

For Theorem \ref{HODGE}, consider the transcendental lattices of $S$ and $X$,
\begin{align*}
T(S) &:= H^{1,1}(S,\Z)^\perp \subset H^2(S,\Z), \\
T(X) &:= H^{2,2}(X,\Z)^\perp \subset H^4(X,\Z),
\end{align*}
as polarised Hodge structures of weights 2 and 4 respectively.

\begin{lem} \label{hodge_conj_lemma}
Let $S$ be a projective K3 surface and $X$ a cubic fourfold.  If a Hodge class $Z \in H^6(S \times X,\Q)$ induces the zero map $T(S)_\Q(-1) \to T(X)_\Q$ then $Z$ is algebraic.
\end{lem}
\begin{proof}
Throughout we use rational cohomology. Classes in $H^0(S) \otimes H^6(X)$ and $H^4(S) \otimes H^2(X)$ are clearly algebraic, so we may assume that $Z$ lies in the $H^2(S) \otimes H^4(X)$ component of the K\"unneth decomposition of $H^6(S\times X)$.

Since $Z$ is a Hodge class, its class in $\Hom_\Q(H^2(S)^*, H^4(X))$ is a morphism of Hodge structures \cite[Lem.~11.41]{voisin_book}. Via the decompositions
\begin{align*}
H^2(S) &= H^{1,1}(S,\Q) \oplus T(S)_\Q, &
H^4(X) &= H^{2,2}(S,\Q) \oplus T(X)_\Q,
\end{align*}
the class $Z\in\Hom^{}_\text{Hdg}(H^2(S)(-1),\,H^4(X))$
therefore defines an element of
\begin{align*}
\Hom^{}_\text{Hdg}(H^{1,1}(S,\Q)(-1),\,H^{2,2}(X,\Q)) &\oplus
\Hom^{}_\text{Hdg}(T(S)_\Q(-1),\,H^{2,2}(X,\Q)) \\
\oplus \Hom^{}_\text{Hdg}(H^{1,1}(S,\Q)(-1),\,T(X)_\Q)
&\oplus \Hom^{}_\text{Hdg}(T(S)_\Q(-1),\,T(X)_\Q).
\end{align*}
It is a standard fact that the second summand vanishes: any morphism $T(S)_\Q(-1) \to H^{2,2}(X,\Q)$ kills $H^{2,0}(S) \subset T(S)_\C$ so has nontrivial kernel. But $T(S)_\Q$ is irreducible \cite[Prop.~2.3]{huybrechts_notes}, so this kernel is all of $T(S)_\Q$.

A similar argument with $H^{3,1}(X)$ shows the third summand vanishes once we know that $T(X)_\Q$ is also irreducible. The proof in \cite[Prop.~2.3]{huybrechts_notes} for $T(S)_\Q$ relies only on the non-degeneracy of the pairing on $H^{1,1}(S,\Q)$ (so that $H^{1,1}(S,\Q) \cap T(S)_\Q = 0$). The pairing on $H^{2,2}(X,\Q)$ is also non-degenerate by the Hodge-Riemann bilinear relations, so the same proof applies.

By hypothesis the component of $Z$ in the fourth summand is zero. So $Z$ lies in the first summand, and is thus algebraic because the Hodge conjecture holds for $S$ and $X$, the latter by \cite{zucker}.  
\end{proof}

The same argument shows that if $\rank(T(S)) \ne \rank(T(X))$ then all Hodge classes on $S\times X$ are products of Hodge classes on $S$ and $X$, so the Hodge conjecture is trivially true for $S\times X$.

\begin{proof}[Proof of Theorem \ref{HODGE}]
Let $\phi_Z\colon T(S)(-1) \to T(X)$ be the Hodge isometry induced by $Z$.  As in the proof of Proposition \ref{mukai-orlov} we can extend this to a Hodge isometry $\Ktop(S) \xrightarrow\sim \Ktop(\A_X)$, so $\Knum(S) \cong \Knum(\A_X)$.  Since $\Knum(S)$ contains a copy of $U$, by Theorem \ref{day} we have $X \in \cC_d$ for some $d$ satisfying \eqref{numerical_condition}.  Thus by Theorem \ref{algcycle} there is a K3 surface $S'$ and an algebraic class $Z' \in H^6(S' \times X, \Q)$ inducing a Hodge isometry $\phi_{Z'}\colon T(S')(-1) \xrightarrow\sim T(X)$.  By \cite{mukai_tata} there is an algebraic class $Z'' \in H^4(S \times S', \Q)$ inducing $\phi_Z \circ \phi_{Z'}^{-1}: T(S') \xrightarrow\sim T(S)$.  Then $Z$ and
\[ Z' \circ Z'' := \pi_{SX*}(\pi_{SS'}^* Z'' \cup \pi_{S'X}^* Z') \]
induce the same map $T(S)(-1) \to T(X)$. Therefore $Z - (Z' \circ Z'')$ is algebraic by Lemma \ref{hodge_conj_lemma}, so $Z$ is too.
\end{proof}